\documentclass[ijoc,leqno,nonblindrev]{informs3hack} 

\SingleSpacedXI

\usepackage{amsfonts}
\usepackage{booktabs}
\usepackage{siunitx}
\usepackage{multirow}
\usepackage{array}
\newcolumntype{C}[1]{>{\centering\arraybackslash}p{#1}}

\usepackage{natbib}
 \bibpunct[, ]{(}{)}{,}{a}{}{,}%
 \def\bibfont{\small}%

\usepackage{placeins}
\usepackage[hypertexnames=false,colorlinks=true,breaklinks=true,bookmarks=true,urlcolor=blue,citecolor=blue,linkcolor=blue,bookmarksopen=false,draft=false]{hyperref}

\TheoremsNumberedThrough     

\EquationsNumberedThrough    

\makeatletter
\let\c@table\c@figure 
\let\ftype@table\ftype@figure 
\makeatother

\DeclareMathOperator{\Diag}{Diag}
\DeclareMathOperator{\diag}{diag}
\DeclareMathOperator{\ldet}{ldet}

\DeclareMathOperator{\Trace}{tr}

\newtheorem{thm}{Theorem}[section]{\bf}{\rm}
\newtheorem{cor}[thm]{Corollary}{\bf}{\rm}
\newtheorem{lem}[thm]{Lemma}{\bf}{\rm}
{\bf}{\rm}
{\bf}{\rm}
\newtheorem{rem}[thm]{Remark}{\bf}{\rm}
{\bf}{\rm}

\begin{document}


\RUNAUTHOR{Chen, Fampa and Lee}

\RUNTITLE{Computing with convex relaxations for MESP}

\TITLE{On computing with some convex relaxations\\ for the maximum-entropy sampling problem}

\ARTICLEAUTHORS{%
\AUTHOR{Zhongzhu Chen}
\AFF{University of Michigan, \EMAIL{zhongzhc@umich.edu}}
\AUTHOR{Marcia Fampa}
\AFF{Universidade Federal do Rio de Janeiro, \EMAIL{fampa@cos.ufrj.br}}
\AUTHOR{Jon Lee}
\AFF{University of Michigan, \EMAIL{jonxlee@umich.edu}}
} 

\ABSTRACT{%
Based on a factorization of an input covariance matrix, we define a mild generalization of an upper bound of Nikolov (2015) and Li and Xie (2020) for the NP-Hard constrained maximum-entropy sampling problem (CMESP). We demonstrate that this
factorization bound is
invariant under scaling and also independent of the particular factorization chosen.
We give a variable-fixing methodology that could be used in a branch-and-bound scheme
based on the factorization bound for
exact solution of CMESP, and we
demonstrate that its ability to fix is independent of the
factorization chosen.
We report on successful experiments with a commercial nonlinear-programming solver.
We further demonstrate that the known ``mixing'' technique
can be successfully used to
combine the factorization bound
with the factorization bound of the complementary CMESP, and also with
the ``linx bound'' of Anstreicher (2020).
}%



\maketitle

%



\section{Introduction}\label{sec:int}

Let $C$ be a symmetric positive semidefinite matrix with rows/columns
indexed from $N:=\{1,2,\ldots,n\}$, with $n >1$.
For $0< s < n$,
we define the \emph{maximum-entropy sampling problem}
\begin{equation}\tag{MESP}\label{MESP}
\begin{array}{ll}
z(C,s):=&\max \left\{\ldet C[S(x),S(x)] ~:~ \mathbf{e}^\top x =s,~ x\in\{0,1\}^n\right\},
\end{array}
\end{equation}
where $S(x)$ denotes the support of $x\in\{0,1\}^n$,
 $C[S,S]$ denotes the principal submatrix indexed by $S$, and
$\ldet$ denotes the natural logarithm of the determinant.

This problem, which has application to contracting environmental monitoring networks (see \cite*{ZidekSunLe2000,CaseltonZidek1984,CaseltonKanZidek1991,Wu-Zidek1992,LeeEnv}, for example),
is a mathematical-programming formulation of the NP-Hard
problem of finding a subset of $s$ random variables from a
Gaussian random $n$-vector having covariance matrix $C$, so
as to maximize the ``information'' (measured by ``differential entropy'');
see \cite*{SW,KLQ,SebWynn,LeeEnv,Kurt_linx}, and the many references they contain.

In practical applications, there are often side constraints (e.g., based on budgetary restriction,
geographical considerations, etc). In this spirit,
for  $A\in\mathbb{R}^{m\times n}$ and $b\in\mathbb{R}^m$,
we further define the \emph{constrained maximum-entropy sampling problem}
\begin{equation}\tag{CMESP}\label{CMESP}
\begin{array}{ll}
z(C,s,A,b):=&\max \left\{\textstyle \vphantom{\sum_{j\in S}} \ldet C[S(x),S(x)] ~:~ \mathbf{e}^\top x =s,~Ax\leq b,~ x\in\{0,1\}^n\right\}.\\
\end{array}
\end{equation}

The main approach for solving \ref{MESP} and \ref{CMESP} to optimality is branch-and-bound.
Lower bounds are calculated by local search, rounding, etc.
Branching is done by deleting a symmetric row/column pair from $C$ for a ``down branch'' and
calculating a Schur complement for an ``up branch''.
Upper bounds are calculated in a variety of ways.
One class of methods uses spectral information; see
\cite*{KLQ,LeeConstrained,HLW,AnstreicherLee_Masked,BurerLee}.
Another class of methods is based on convex relaxations; see
\cite*{AFLW_Using,Anstreicher_BQP_entropy,Kurt_linx}.
Recently, \cite*{Nikolov} and \cite*{Weijun} developed a convex relaxation based on
spectral information.

It is easy to check that
\[
z(C,s,A,b)=z(C^{-1},n-s,-A,b-A\mathbf{e}) + \ldet C~.
\]
\noindent 
So we have a notion of a
\emph{complementary} \ref{CMESP} problem
\begin{align*}\tag{CMESP-comp}\label{CMESP-comp}
&z(C^{-1},n-s,-A,b-A\mathbf{e})= \\
&\qquad  \max \left\{\textstyle \vphantom{\sum_{j\in S}} \ldet C^{-1}[S(x),S(x)] ~:~ \mathbf{e}^\top x =n-s,~A(\mathbf{e}-x)\leq b,~ x\in\{0,1\}^n\right\},
\end{align*}
 and \emph{complementary} bounds
 (i.e., bounds for the
 complementary problem plus $\ldet C$)
 immediately give us bounds on
 $z(C,s,A,b)$.
 Some upper bounds on $z(C,s,A,b)$ also shift by  $\ldet C$  under complementing,
 in which case there is no additional value in computing the
 complementary bound.

 It is also easy to check that
\[
z(C,s,A,b)=z(\gamma C,s,A,b) - s \log \gamma,
\]
\noindent where the \emph{scale factor} $\gamma>0$.
So we have a notion of a
\emph{scaled} \ref{CMESP} problem defined by the data
 $\gamma C$, $s$, $A$, $b$, and \emph{scaled} bounds
 (i.e., bounds for the
 scaled problem minus $s \log \gamma$)
 immediately give us bounds on
 $z(C,s,A,b)$.
 Some upper bounds on $z(C,s,A,b)$ also shift by  $-s \log \gamma$  under scaling,
 in which case there is no additional value in computing the
scaled bound.

In Section \ref{sec:bounds}, we discuss some upper bounds for \ref{CMESP}.
In particular, we present the ``factorization bound'' for \ref{CMESP},
as a convex formulation \ref{DFact} of the Lagrangian dual of a nonconvex formulation \ref{Fact},
and some of its important properties from a computational point of view.
In particular:
\begin{itemize}
\item We demonstrate that the factorization bound changes by the same amount as $z(C,s,A,b)$, when $C$ is scaled by some $\gamma>0$.
\item  We give a variable fixing methodology based on a feasible solution
of  \ref{DFact}. Variable fixing has also been studied for \ref{CMESP} in the context of other convex relaxations; see \cite*{AFLW_IPCO,AFLW_Using,AFLW_Remote,Anstreicher_BQP_entropy,Kurt_linx}.
\item  We demonstrate that the factorization bound and its ability to fix variables, based on a matrix factorization $C=FF^\top$, is independent of the factorization (and so mathematically, it provides the same bound as that of
\cite*{Nikolov} and \cite*{Weijun}).
\item  We demonstrate that the factorization bound dominates the
well-known ``spectral bound''.
\item  Although it is possible to directly attack \ref{DFact}
to calculate the factorization bound, it is much more fruitful to
work with  \ref{DDFact}, a convex formulation of the Lagrangian dual of
\ref{DFact}. In connection with this, we provide a mechanism to get a minimum-gap feasible solution
of \ref{DFact}, relative to a  feasible (and possibly non-optimal)
solution of   \ref{DDFact} (which is useful for getting a true upper bound for \ref{CMESP}). We also describe how to get the gradient of the objective function of \ref{DDFact} (under a technical condition), which is necessary for
applying any reasonable technique for efficiently solving  \ref{DDFact}.
\item We review the \ref{linx} bound for \ref{CMESP},
and we present some of its key properties that are useful for computing.
\item We review  the ``mixing'' bound of \cite*{Mixing}, and we work out a dual formulation for it, as well as a fixing methodology in an important case that generalizes what we can do for \ref{DDFact},
complementary \ref{DDFact}, and \ref{linx}.
\end{itemize}

In Section \ref{sec:impl},
we discuss the numerical experiments where, using a commercial nonlinear-programming solver, we calculate  upper bounds for benchmark instances of \ref{MESP} from the literature with the three relaxations presented, namely \ref{DDFact}, complementary \ref{DDFact} and \ref{linx}, and  with the ``mixing'' strategy described in \cite*{Mixing}. Generally, we found that a
commercial nonlinear-programming solver is quite viable for
our relaxations, even for \ref{DDFact} which may have nondifferentiabilty. Our main findings:
\begin{itemize}
     \item We compared integrality gaps given by the difference between the upper bounds computed with the relaxations and
     the lower bounds computed with a greedy/interchange heuristic
     or with the optimal value when we could obtain it by branch-and-bound. We found  that all the three relaxations,  \ref{DDFact}, complementary \ref{DDFact}, and \ref{linx}, achieve the best bounds for some of the instances, however \ref{DDFact} and \ref{linx} achieve together most of the best bounds.
    \item We compared the times to solve the relaxations and analyzed the impact of two factors in the times: the smoothness of the objective functions of the relaxations and the ranks of the covariance matrices. The possibility of  \ref{DDFact} and complementary \ref{DDFact} encountering points at which the objective function is  nondifferentiable
    led us to the application of a BFGS-based algorithm to solve these relaxations, and the smoothness of \ref{linx} results in the application of a  Newton-based algorithm to solve it. We see a better convergence of the Newton-based algorithm, resulting in best times for \ref{linx} and with less variability among the different values of $s$, except for our largest instance with $n=2000$ and a covariance matrix with rank $r=949$. In this case, \ref{linx} presents the drawback of dealing with an order-$n$ matrix in its objective function, while \ref{DDFact} deals with an order-$r$ matrix.
    \item We demonstrated how the ``mixing'' procedure can decrease the bounds obtained with the three relaxations when we mix two relaxations that obtain very similar bounds when applied separately. This is mostly observed when we mix \ref{DDFact} and \ref{linx}. For the majority of the instances, \ref{DDFact} presents better bounds for values of $s$ up to an intermediate value, and for larger values of $s$, \ref{linx} presents better bounds. For values of $s$ close to this intermediate value, the mixing strategy effectively decreases the bound obtained by each relaxation.
    \item We demonstrated how the fixing methodology can fix a significant number of variables, especially when applied iteratively, to our largest instance.
\end{itemize}

In Section \ref{sec:discussion}, we summarize our results
and point to future work.
\medskip

\noindent {\bf Notation.}
 We let $\mathbb{S}^n_+$ (resp., $\mathbb{S}^n_{++}$) denote the set of
 positive semidefinite (resp., definite) symmetric matrices of order $n$.
We denote by $\lambda_\ell(M)$  the $\ell$-th greatest eigenvalue of
a matrix $M\in \mathbb{S}^n_+$~. We denote an all-ones  vector
by $\mathbf{e}$ and the $i$-th standard unit vector by $\mathbf{e}_i$~.
For matrices $A$ and $B$ with the same shape,
$A\circ B $ is the Hadamard (i.e., element-wise) product,
and $A\bullet B:=\Trace(A^\top B)$ is the matrix dot-product.
For a matrix $A$, we denote row $i$ by $A_{i\cdot}$ and
column $j$ by $A_{\cdot j}$~.

\section{Upper bounds}\label{sec:bounds}

There are a wide variety of upper bounding methods for \ref{CMESP}.
The tightest bounds, from a practical computational viewpoint,
seem to be  Anstreicher's ``linx bound'' and the ``factorization bound''.
In this section, we describe and develop these bounds, with an
eye on practical and efficient computation.

\subsection{Fact}
We begin by introducing a mild generalization of a
nonconvex programming bound developed by
 \cite*{Nikolov} and \cite*{Weijun}.
Suppose that the rank of $C$ is $r\geq s$.
Then we factorize $C=FF^\top$,
with $F\in \mathbb{R}^{n\times k}$, for some $k$ satisfying $r\le k \le n$.
We note that this could be a
Cholesky-type factorization (i.e., $k:=r$ and $F$ lower triangular), as in  \cite*{Nikolov} and \cite*{Weijun}.
But it could alternatively be derived from
a spectral decomposition of $C$; that is,
$C=\sum_{\ell=1}^r \lambda_\ell v_\ell v_\ell^\top$,
where we  put $\sqrt{\lambda_\ell}v_\ell$ as column $\ell$ of $F$,
$\ell=1,\ldots,k:=r$. Another very useful possibility is to
let $k:=n$, and choose $F$ to be the matrix square root, $C^{1/2}$, which is
always symmetric.

Next, for $x\in [0,1]^n$, we define
$
   F(x) := \sum_{j\in N} F_{j\cdot}^\top F_{j\cdot} ~x_j~ =F^\top \Diag(x)F
$
and
	\[
\begin{array}{ll}
z_{{\tiny\mbox{Fact}}}(C,s,A,b;F):=&\max \displaystyle \sum_{\ell=1}^s \log\left( \lambda_{\ell}  (F(x))\right)\\
  &\mbox{subject to:}\\
     &\quad\mathbf{e}^{\top}x=s,~ Ax\leq b,~ 0\leq x\leq \mathbf{e}.
\end{array}
\tag{Fact}\label{Fact}
\]

It is easy to check that $z(C,s,A,b) \leq z_{{\tiny\mbox{Fact}}}(C,s,A,b;F)$, for any factorization of $C$ (c.f. \cite*{KLQ}).
Unfortunately,
\ref{Fact} is not a convex program, so it is not a practical relaxation to work with.

\subsection{DFact}\label{sec:DFact}

We define
\[
f(\Theta,\nu, \pi,\tau):= -\sum_{\ell=k-s+1}^k \log\left(\lambda_{\ell} \left(\Theta\right)\right)
+ \nu^\top \mathbf{e}  + \pi^\top b +\tau s - s,
\]
and the \emph{factorization bound}
\[
\begin{array}{ll}
	z_{{\tiny\mbox{DFact}}}(C,s,A,b;F):=& 	\min~ f(\Theta,\nu, \pi,\tau)\\
     & \mbox{subject to:}\\
     &\quad \diag(F \Theta F^\top) + \upsilon - \nu  - A^\top \pi - \tau\mathbf{e}=0,\\
&\quad \Theta\succ 0, ~\upsilon\geq 0, ~\nu\geq 0, ~\pi\geq 0.
\end{array}
\tag{DFact}\label{DFact}
\]

\ref{DFact} is equivalent to the Lagrangian dual of \ref{Fact}, and it is a convex program.
The objective function of \ref{DFact}  is analytic at
every point $(\hat\Theta,\hat\upsilon,\hat\nu,\hat\pi,\hat\tau)$ for
which $\lambda_{k-s}(\hat\Theta) > \lambda_{k-s+1}(\hat\Theta)$.
In fact, we have seen in our experiments, a good solver can
get to an optimum, even when this condition fails.

It turns out that the factorization bound for \ref{MESP} has a close relationship with the spectral bound of \cite*{KLQ}:
$\sum_{\ell=1}^s \log \lambda_\ell(C)$.
First, we establish that like the spectral bound for \ref{MESP}, the
factorization bound for \ref{CMESP} is invariant under multiplication of $C$ by a scale factor $\gamma$, up to the
additive constant $-s\log \gamma$, a property that is
\emph{not} shared with other convex-optimization bounds.

\begin{thm}\label{thm:scale}
For all $\gamma>0$ and factorizations $C=FF^\top$,
we have
\[
z_{{\tiny\mbox{DFact}}}(C,s,A,b;F) =  z_{{\tiny\mbox{DFact}}}(\gamma C,s,A,b; \sqrt{\gamma}F)  -s \log \gamma .
\]
\end{thm}

\proof{Proof.}
We simply observe that for every feasible solution  $(\hat{\Theta},\hat{\upsilon},\hat{\nu},\hat{\pi},\hat{\tau})$
of \ref{DFact},
we have that $(\frac{1}{\gamma}\hat{\Theta},\hat{\upsilon},\hat{\nu},\hat{\pi},\hat{\tau})$ is a feasible solution
of \ref{DFact} with $F$ replaced by $\sqrt{\gamma}F$.
Then we observe that  $\lambda_{\ell}\left(\frac{1}{\gamma}\hat{\Theta}\right) = \frac{1}{\gamma}\lambda_{\ell}(\hat{\Theta})$, for all $\ell$.
This mapping between feasible solutions is a bijection, so the result follows.
\Halmos\endproof

Because the factorization bound shifts by the same amount as $z(C,s,A,b)$,
under scaling of $C$ by $\gamma$,
we cannot improve on the factorization bound by scaling.
In contrast, the linx bound is very sensitive to the choice of the scale factor,
and while we can compute an optimal scale factor for the linx bound (see \cite*{Mixing}),
it is a significant computational burden to do so.

Next, we present another useful result that guides practical usage.

\begin{thm}\label{thm:FactDoesntMatter}
Let $C=F_jF_j^\top$, for $j=1,2$, be two different factorizations of $C$, and let  $(\hat\Theta_1,\hat\upsilon,\hat\nu,\hat\pi,\hat\tau)$ be a feasible solution to \ref{DFact}, for $F:=F_1$. Then, there is a feasible solution $(\hat\Theta_2,\hat\upsilon,\hat\nu,\hat\pi,\hat\tau)$ to \ref{DFact}, for $F:=F_2$, such that $f(\hat\Theta_1,\hat\nu,\hat\pi,\hat\tau)=f(\hat\Theta_2,\hat\nu,\hat\pi,\hat\tau)$.
\end{thm}

\proof{Proof.}%
Let $r$ be the rank of $C$, and let $C=\sum_{\ell=1}^n \lambda_\ell u_\ell u_\ell^\top$ be a spectral
decomposition of $C$. Suppose that $C=FF^\top$, with $F\in\mathbb{R}^{n\times k}$, and $r\leq k\leq n$. Our preliminary goal is to build a special singular-value decomposition of $F$.

Let $\sigma_\ell:=\sqrt{\lambda_\ell}$, for $1\leq \ell \leq k$.
Now define $v_\ell\in\mathbb{R}^k$, for $1\leq \ell\leq r$ by
$v_\ell:= \frac{1}{\sigma_\ell} F^\top u_\ell$.

We can easily check that for $1\leq i \leq \ell \leq r$, we have
\[
v_i^\top v_\ell = \frac{1}{\sigma_i \sigma_\ell}u_i^\top F F^\top u_\ell = \frac{1}{\sigma_i \sigma_\ell}u_i C u_\ell
= \frac{\lambda_\ell}{\sigma_i \sigma_\ell} u_i^\top u_\ell =
\left\{
  \begin{array}{ll}
  1 , & \hbox{ for $i=\ell$;} \\
  0 , & \hbox{ for $i<\ell$.}
  \end{array}
\right.
\]
That is,
$\{ v_\ell ~:~ 1\leq \ell \leq r\}$ is a set of $r$ orthonormal vectors in
$\mathbb{R}^k$. So, for $r<\ell\leq k$, we can now choose $v_i$ so as to complete $\{ v_\ell ~:~ 1\leq \ell \leq r\}$ to an orthonormal basis of $\mathbb{R}^k$.

Next, we have
\[
\sum_{\ell=1}^k \sigma_\ell u_\ell v_\ell^\top =
\sum_{\ell=1}^r \sigma_\ell u_\ell v_\ell^\top =
\sum_{\ell=1}^r  u_\ell u_\ell^\top F =
\sum_{\ell=1}^n  u_\ell u_\ell^\top F =
I_n F =
F
\]
(note that above we have use the fact that $\| u_\ell^\top F \|^2 =  u_\ell^\top F F^\top u_\ell = \lambda_\ell=0$ for $r<\ell \leq n$, which implies that $u_\ell^\top F=0^\top_k$ for  $r<\ell\leq n$),
and so we can conclude that
$F=\sum_{\ell=1}^k \sigma_\ell u_\ell v_\ell^\top$ is a singular-value decomposition for $F$.

The important takeaway is that the $u_\ell\in\mathbb{R}^n$ ($1\leq \ell \leq n$) and the nonzero $\sigma_\ell$
 ($1\leq \ell \leq r$) in the singular-value decomposition that we constructed for $F$ only depend on $C$, not on the
 particular factorization $C=FF^\top$.

It is convenient now to establish that in a factorization matrix $F$, we
can without loss of generality take $k=n$, by appending 0 columns to $F$ if needed, and this will not affect the
bound $z_{{\tiny\mbox{DFact}}}(C,s,A,b;F)$.
Let $\bar{F}:= [F ~\vert~ 0_{n\times(n-k)}]$,
and consider
\[
\bar{\Theta} =
\left(
  \begin{array}{cc}
\hat \Theta & \times  \\
\times  & \times \\
  \end{array}
\right)\in  \mathbb{S}^n_+~.
\]
It is trivial to check that $\bar F \bar \Theta \bar F^\top = F \hat \Theta F^\top$.
And by Cauchy's eigenvalue interlacing inequalities (see \cite*{HJBook}, for example), we have $\lambda_{\ell+n-k}(\bar \Theta) \leq \lambda_\ell (\hat\Theta)$,
for $1\leq \ell \leq k$.
Therefore, we have $z_{{\tiny\mbox{DFact}}}(C,s,A,b;\bar F) \geq
z_{{\tiny\mbox{DFact}}}(C,s,A,b;F)$. In the other direction,
suppose that $\hat \Theta\in \mathbb{S}^k_+$~.
Now define
\[
\bar{\Theta} :=
\left(
  \begin{array}{cc}
\hat \Theta & 0^\top \\
0 & \lambda_1(\hat\Theta) I_{n-k} \\
  \end{array}
\right)\in  \mathbb{S}^n_+~.
\]
As above, we have  $\bar F \bar \Theta \bar F^\top = F \hat \Theta F^\top$.
And by construction, we have
$\lambda_{\ell+n-k}(\bar \Theta) = \lambda_\ell (\hat \Theta)$,
for $1\leq \ell \leq k$. And therefore, we have $z_{{\tiny\mbox{DFact}}}(C,s,A,b;\bar{F}) \leq
z_{{\tiny\mbox{DFact}}}(C,s,A,b;F)$.

With this we can now conclude that if we have two  different factorizations of $C$,
say $C=F_j F_j^\top$ for $j=1,2$, we can without loss of generality assume
for each that $F_j$ has $k=n$ columns, and further that
we can choose singular-value decompositions of the form
$F_j=U \Sigma V_j^\top$~, where here we now take $U$, $\Sigma$, $V_1$ and $V_2$
to all be $n\times n$.

Now, right multiplying $U \Sigma V_2^\top=F_2$ by
$V_2 V_1^\top$, we get $U\Sigma V_2^\top V_2 V_1^\top= F_2 V_2 V_1^\top$, and so we have
$F_1=F_2 V_2 V_1^\top$.

Finally, for $\Theta_1\succ 0$, we have
$
F_1 \Theta_1 F_1^\top = F_2 V_2 V_1^\top \Theta_1 V_1 V_2^\top F_2^\top~,
$
and so by taking $\Theta_2:= V_2 V_1^\top \Theta_1 V_1 V_2^\top$, we get
$F_1 \Theta_1 F_1^\top= F_2 \Theta_2 F_2^\top$, with $\Theta_2$ being similar to $\Theta_1$.
Therefore, we have that $\lambda_\ell(\Theta_1)=\lambda_\ell(\Theta_2)$, for all $\ell$,
and so
we can transform any feasible solution of \ref{DFact} with respect to factor $F:=F_1$
into a  feasible solution of \ref{DFact} having the same objective value
with respect to factor $F:=F_2$. The result follows.
\Halmos\endproof

\begin{cor}\label{cor:FactDoesntMatter}
The value of the factorization bound is independent of the particular factorization.
\end{cor}

\begin{rem}
 Corollary \ref{cor:FactDoesntMatter} follows directly from Theorem \ref{thm:FactDoesntMatter}. We note that the proof of Theorem \ref{thm:FactDoesntMatter} not only confirms the statement in Corollary \ref{cor:FactDoesntMatter}, but also presents  a methodology for constructing a feasible solution $(\hat\Theta_2,\hat\upsilon,\hat\nu,\hat\pi,\hat\tau)$ to \ref{DFact} for a given factorization of $C$, from a feasible solution  $(\hat\Theta_1,\hat\upsilon,\hat\nu,\hat\pi,\hat\tau)$ to \ref{DFact} for any other factorization, where both solutions have the same objective value with respect to the corresponding factor. In Section \ref{sec:DDFact}, we also present a short proof for Corollary \ref{cor:FactDoesntMatter}.
\end{rem}

Next, we establish that the factorization bound for \ref{MESP} dominates the spectral bound  for \ref{MESP}. While the spectral bound is
much cheaper to compute, because of this result, there is never any point
of computing the spectral bound if we have already computed the factorization bound.
In another way of thinking, if the spectral bound comes close to
allowing us to discard a  subproblem in the context of branch-and bound, it
should be well worth computing the factorization bound to attempt to
 discard the subproblem.



\begin{thm}
Let $C\in\mathbb{S}^n_{+}$, with $r:=\mbox{rank}(C)$, and $s\leq r$. Then,
for all factorizations $C=F F^\top$, we have
\[
z_{{\tiny\mbox{DFact}}}(C,s,\cdot,\cdot\,\, ;F)\leq  \sum_{\ell=1}^s  \log \lambda_\ell(C)~.
\]
\end{thm}

\proof{Proof.}

Let $C=\sum_{\ell=1}^r \lambda_\ell(C) u_\ell u_\ell^\top$ be a spectral
decomposition of $C$.
By Theorem \ref{thm:FactDoesntMatter},
it suffices to take $F$ to be the
symmetric matrix $\sum_{\ell=1}^r \sqrt{\lambda_\ell(C) } u_\ell u_\ell^\top$.

We consider the solution for \ref{DFact} given by:   $\hat{\Theta}:=C^{\dagger} + \frac{1}{\lambda_r(C)}\left(I-CC^{\dagger}\right)$, where $C^{\dagger}:=\sum_{\ell=1}^r \frac{1}{\lambda_\ell(C)}  u_\ell u_\ell^\top$ is the  Moore-Penrose pseudoinverse  of $C$,  $\hat{\upsilon}:=\mathbf{e} - \diag(F\hat\Theta F^\top)$,
$\hat{\nu}:=0$,  $\hat{\pi}:=0$, and $\hat{\tau}:=1$.
We can verify that the $r$ least eigenvalues of $\hat\Theta$ are $\frac{1}{\lambda_1(C)},\frac{1}{\lambda_2(C)},\ldots,\frac{1}{\lambda_r(C)}$ and the $n-r$ greatest eigenvalues are all equal to  $\frac{1}{\lambda_r(C)}$.
Therefore, $\hat\Theta$ is positive definite.

The equality constraint of \ref{DFact} is clearly satisfied at this solution. Additionally, we can verify that $F\hat\Theta F^\top = \sum_{\ell=1}^r u_\ell u_\ell^\top$. As the positive semidefinite matrix  $\sum_{\ell=r+1}^n u_\ell u_\ell^\top$ is equal to  $I - \sum_{\ell=1}^r u_\ell u_\ell^\top$, we conclude that  $\diag(F\hat\Theta F^\top)\leq \mathbf{e}$. Therefore, $\hat{\upsilon}\geq 0$, and the solution constructed is a feasible solution to \ref{DFact}. Finally, we can see that the   objective value
of this solution is equal to the spectral bound. The result then follows.
\Halmos\endproof

\begin{rem}
 We note that when $C$ is nonsingular,
then $F=C^{1/2}$, and
using the symmetry of $C^{1/2}$, it is easy to directly check that
with  $\hat{\Theta}:=C^{-1}$, $\hat{\tau}:=1$,
$\hat{\upsilon}:=\hat{\nu}:=0$, and  $\hat{\pi}:=0$, we have a feasible
solution of \ref{DFact} with objective value
equal to the spectral bound.
\end{rem}

Next, we consider variable fixing, in the context of solving \ref{CMESP}.

\begin{thm} \label{thm:fixFact}
 Let
\begin{itemize}
    \item $\mbox{LB}$ be the objective-function value of a feasible solution for \ref{CMESP},
    \item $(\hat{\Theta},\hat{\upsilon},\hat{\nu},\hat{\pi},\hat{\tau})$ be a feasible solution for \ref{DFact} with objective-function value $\hat{\zeta}$.
\end{itemize}
Then, for every optimal solution $x^*$
for \ref{CMESP}, we have:
\[
\begin{array}{ll}
x_j^*=0, ~ \forall ~ j\in N \mbox{ such that } \hat{\zeta}-\mbox{LB} < \hat{\upsilon}_j~,\\
x_j^*=1, ~ \forall ~ j\in N \mbox{ such that } \hat{\zeta}-\mbox{LB} < \hat{\nu}_j~.\\
\end{array}
\]
\end{thm}

\proof{Proof.}
Consider \ref{Fact} with the additional constraint $x_i=1$. The dual becomes then,
\begin{equation}\label{modDfact}
\begin{array}{ll}
 \min& -
\displaystyle \sum_{\ell=k-s+1}^k \log\left(\lambda_{\ell} \left(\Theta\right)\right)
+ \nu^\top \mathbf{e}  + \pi^\top b +\tau s - s- \omega\\
     &\mbox{subject to:}\\
     & \diag(F \Theta F^\top) + \upsilon - \nu  - A^\top \pi - \tau\mathbf{e} + \omega \mathbf{e}_j=0,\\
& \Theta\succ 0, ~\upsilon\geq 0, ~\nu\geq 0, ~\pi\geq 0.
\end{array}
\end{equation}
where $\omega$ is the new dual variable.
Notice that, as long as
 $\hat{\upsilon}_j-\omega \geq 0$,
 $(\hat{\Theta},\hat{\upsilon}-\omega\mathbf{e}_j,\hat{\nu},\hat{\pi},\hat{\tau},\omega)$ is a feasible solution of the modified dual, with objective value
 $\hat{\zeta}-\omega$. So, to minimize the
 objective value of our feasible solution of the modified dual,
 we set $\omega$ equal to $\hat{\upsilon}_j$. We conclude that $\hat{\zeta}-\hat{\upsilon}_j$
 is an upper bound on the objective value of every solution of \ref{CMESP} that satisfies $x_j=1$. So if $\hat{\zeta}-\hat{\upsilon}_j<LB$, then
 no optimal solution of \ref{CMESP}  can have $x_j=1$.

Similarly, consider \ref{Fact} with the additional constraint $x_j=0$. In this case, the new dual problem is equivalent to \eqref{modDfact}, except that the objective function does not have the term $-\omega$. Therefore,  as long as
 $\hat{\nu}_j+\omega \geq 0$,
 $(\hat{\Theta},\hat{\upsilon},\hat{\nu}+\omega\mathbf{e}_j,\hat{\pi},\hat{\tau},\omega)$ is a feasible solution of this modified dual with objective value
 $\hat{\zeta}+\omega$, and to minimize the
 objective value of the feasible solution,
 we set $\omega$ equal to $-\hat{\nu}_j$. Now, we conclude that $\hat{\zeta}-\hat{\nu}_j$
 is an upper bound on the objective value of every solution of \ref{CMESP} that satisfies $x_j=0$. So if $\hat{\zeta}-\hat{\nu}_j<LB$, then
 no optimal solution of \ref{CMESP}  can have $x_j=0$.
\Halmos\endproof

\begin{rem}
We note that Theorem \ref{thm:FactDoesntMatter} implies that all factorizations
$C=FF^\top$ have the same power to fix variables.
\end{rem}

It is quite possible to develop a direct nonlinear-programming algorithm to attack \ref{DFact}.
For any reasonably-fast algorithm, we
would need the gradient of  \ref{DFact}  objective function. Toward this,
we consider a spectral decomposition of $\hat\Theta$, that is $\hat\Theta=\sum_{\ell=1}^{k} \lambda_{\ell}(\hat\Theta) u_{\ell}(\hat\Theta) u_{\ell}(\hat\Theta)^{\top}$.
If $\lambda_{k-s}(\hat\Theta)>\lambda_{k-s+1}(\hat\Theta)$, then, using \cite*[Theorem 3.1]{Tsing},
the gradient of the objective function of \ref{DFact} is given by
\begin{align*}
&\frac{\partial f(\hat\Theta,\hat\nu,\hat\pi,\hat\tau)}{\partial \Theta} =   -\sum_{\ell=k-s+1}^k\frac{1}{\lambda_{\ell}(\hat\Theta)}u_\ell(\hat\Theta) u_\ell(\hat\Theta)^\top;\quad\frac{\partial f(\hat\Theta,\hat\nu,\hat\pi,\hat\tau)}{\partial \nu} = \mathbf{e};\\[5pt]
&\frac{\partial f(\hat\Theta,\hat\nu,\hat\pi,\hat\tau)}{\partial \pi} = b;\quad \frac{\partial f(\hat\Theta,\hat\nu,\hat\pi,\hat\tau)}{\partial \tau} = s.
\end{align*}
Without the technical condition $\lambda_{k-s}(\hat\Theta)>\lambda_{k-s+1}(\hat\Theta)$,
the formulae above still gives a subgradient of $f$
(see \cite*{Nikolov} for details).

\subsection{DDFact}\label{sec:DDFact}

While it turns out that the bound given by \ref{DFact} is generally quite good,
and it has the potential to fix variables at 0/1 values via Theorem \ref{thm:fixFact},
the model \ref{DFact} is not easy to solve directly. We instead present its (equivalent)
Lagrangian dual,
\ref{DDFact}, which is much easier to work with computationally.

\begin{lem}\label{Ni13}(see \cite*[Lemma 13]{Nikolov})
 Let $\lambda\in\mathbb{R}_+^k$ with $\lambda_1\geq \lambda_2\geq \cdots\geq \lambda_k$ and let $0<s\leq k$. There exists a unique integer $\iota$, with $0\leq \iota< s$, such that
 \begin{equation*}\label{reslemiota}
 \lambda_{\iota}>\frac{1}{s-\iota}\sum_{\ell=\iota+1}^k \lambda_{\ell}\geq \lambda_{\iota +1} ,
 \end{equation*}
 with the convention $\lambda_0=+\infty$.
\end{lem}

Suppose that  $\lambda\in\mathbb{R}^k_+$, and assume that
$\lambda_1\geq\lambda_2\geq\cdots\geq\lambda_k$. Given an integer $s$ with $0<s\leq k$,
let $\iota$ be the unique integer defined by Lemma \ref{Ni13}. We define
\begin{equation*}\label{defGamma}
\phi_s(\lambda):=\sum_{\ell=1}^{\iota} \log\left(\lambda_\ell\right) + (s - \iota)\log\left(\frac{1}{s-{\iota}} \sum_{\ell=\iota+1}^{k}
\lambda_\ell\right).
\end{equation*}
Next, for $X\in\mathbb{S}_{+}^k$, we define
$\Gamma_s(X):= \phi_s(\lambda_1(X),\ldots,\lambda_k(X))$.
Finally, we define
\[
\begin{array}{ll}
	z_{{\tiny\mbox{DDFact}}}(C,s,A,b;F):= 	&\max~ \Gamma_s(F(x)) \\
     &\mbox{subject to:}\\
     &\quad \mathbf{e}^\top x=s,~ Ax\leq b,~
 0\leq x\leq \mathbf{e}.
\end{array}
\tag{DDFact}\label{DDFact}
\]
It is a result of \cite*{Nikolov} that \ref{DDFact}
is a convex program, and that it is in fact
equivalent to the Lagrangian dual of \ref{DFact}.
Checking a Slater's condition, we
have that $	z_{{\tiny\mbox{DDFact}}}(C,s,A,b;F)=	z_{{\tiny\mbox{DFact}}}(C,s,A,b;F)$.
The advantage of solving \ref{DDFact} instead of \ref{DFact} is that it has
many fewer variables. But, variable fixing (see Theorem \ref{thm:fixFact})
relies on a good feasible solution of \ref{DFact}. Moreover,
certifying the quality of a feasible solution of \ref{DDFact} also requires
a good feasible solution of \ref{DFact}.
Motivated by these points, we show how to construct a feasible solution of \ref{DFact} from a feasible solution $\hat x$ of \ref{DDFact}  with finite objective value, with the goal of producing a small gap.

We  consider the spectral
decomposition $F(\hat{x})=\sum_{\ell=1}^{k} \hat \lambda_\ell \hat u_\ell \hat u_\ell^\top,$
with $\hat \lambda_1\geq\hat \lambda_2\geq\cdots\geq \hat \lambda_{\hat r}>\hat \lambda_{\hat{r}+1}=\cdots=\hat \lambda_k=0$.  Notice that $\mbox{rank}(F(\hat x))= \hat{r}\geq s$.
Following \cite*{Nikolov}, we  define
$\hat{\Theta}:=\sum_{\ell=1}^{k} {\hat \beta}_\ell \hat{u}_\ell \hat{u}_\ell^\top$,
where
\begin{equation}\label{defbeta}
\hat{\beta}_\ell:=\left\{
\begin{array}{ll}
        \textstyle 1/\hat{\lambda}_\ell,
       &~1\leq \ell\leq \hat{\iota};\\
     1/\hat{\delta},&~\hat{\iota}<\ell\leq \hat{r};\\
     (1+\epsilon)/\hat{\delta},&~\hat{r}<\ell\leq k,
\end{array}\right.
\end{equation}
 for any $\epsilon>0$, where $\hat{\iota}$ is the unique integer defined  in Lemma \ref{Ni13} for $\lambda_\ell=\hat{\lambda}_\ell$, and
$
\hat \delta:=\frac{1}{s-\hat \iota}\sum_{\ell=\hat \iota+1}^{k}\hat \lambda_\ell
$~.
From Lemma \ref{Ni13}, we have that $\hat\iota<s$.  Then,

\begin{equation}\label{sumbeta}
- \sum_{\ell=1}^s \log\left(\hat{\beta}_{\ell}\right)=  \sum_{\ell=1}^{\hat{\iota}} \log\left(\hat{\lambda}_{\ell}\right) + (s-\hat{\iota})\log(\hat{\delta})= \Gamma_s(F(\hat{x})).
\end{equation}

The minimum duality gap between $\hat x$ in \ref{DDFact} and feasible solutions
of \ref{DFact} of the form $(\hat\Theta,\upsilon,\nu,\pi,\tau)$
is the optimal value of
\begin{equation}\label{mingapprob}\tag{$G(\hat\Theta)$}
\begin{array}{l}
\min~
 \nu^\top \mathbf{e}  + \pi^\top b +\tau s - s\\
      \mbox{subject to:}\\
     \quad  \upsilon - \nu  - A^\top \pi - \tau\mathbf{e}= - \diag(F \hat \Theta F^\top) ,\\
\quad \upsilon\geq 0, ~\nu\geq 0, ~\pi\geq 0.
\end{array}
\end{equation}
Note that \ref{mingapprob} is always feasible (e.g.,
$\upsilon:=0$, $\nu:=\diag(F \hat \Theta F^\top)$, $\pi:=0$, $\tau:=0$). Also,  \ref{mingapprob} has a
simple closed-form solution for \ref{MESP}, that is when there are no $Ax\leq b$  constraints (see \cite*{Weijun}).

Next, we restrict our attention to \ref{MESP}, and we consider the behavior of
the optimal value of \ref{mingapprob} as a function of $\epsilon$.
Let $x^*\in\{0,1\}^n$ be the support vector of the $s$ greatest elements of $\diag(F\hat \Theta F^\top)$.
Then the optimal value of \ref{mingapprob} (and its dual) is $\diag(F\hat \Theta F^\top)^\top x^*-s =\sum_{\ell=1}^k \diag( (F\hat U)^\top \Diag(x^*) (F \hat U))_\ell \hat\beta_\ell -s$, where $\hat U$ is the matrix
having $\ell$-th column equal to $\hat u_\ell$, for $1\leq \ell\leq k$. It is easy to see
that the diagonal elements of $(F\hat U)^\top \Diag(x^*) (F \hat U)$ are nonnegative.
Therefore, with $x^*$ fixed, $\sum_{\ell=1}^k \diag( (F\hat U)^\top \Diag(x^*) (F \hat U))_\ell \hat\beta_\ell$
is non-decreasing in $\epsilon$. Now the optimal value of the dual of \ref{mingapprob}, is the
point-wise $\max$, over the choices of $x^*\in\{0,1\}^n$ satisfying $\mathbf{e}^\top x^*=s$.
So, the optimal value of the dual of \ref{mingapprob} is the point-wise max of linear functions, each of which is
non-decreasing in  $\epsilon$. And so the  optimal value of the dual is also non-decreasing in $\epsilon$.

For developing a reasonable nonlinear-programming algorithm
for \ref{DDFact}, we need an expression for the gradient of its
 objective function.

\begin{thm}
Let $F(\hat x)=\sum_{\ell=1}^{k} \hat\lambda_{\ell} \hat u_{\ell} \hat u_{\ell}^{\top}$ be a spectral decomposition of $F(\hat x)$.
Let $\hat \iota$ be the value of
$\iota$ in Lemma \ref{Ni13},
where $\lambda$ in Lemma \ref{Ni13}
is $\hat{\lambda}:=\lambda(F(\hat x))$.
 If $\frac{1}{s-\hat \iota}\sum_{\ell=\hat \iota+1}^k \hat \lambda_{\ell}> \hat \lambda_{\hat\iota +1}$, then, for
$j=1,2,\ldots,n$,
\begin{equation*}
      \frac{\partial}{\partial x_j} \Gamma_s(F(\hat x))
    =  \sum_{\ell=1}^{\hat\iota} \frac{1}{\hat\lambda_\ell}(F_{j\cdot} \hat u_\ell)^2
    + \sum_{\ell=\hat\iota+1}^{k} \frac{s-\hat\iota}{\sum_{i=\hat\iota+1}^{k}\hat\lambda_i}  (F_{j\cdot} \hat u_\ell)^2~.
\end{equation*}
\end{thm}

\proof{Proof.}
Under the hypothesis $\frac{1}{s-\hat \iota}\sum_{\ell=\hat \iota+1}^k \hat \lambda_{\ell}> \hat \lambda_{\hat\iota +1}$,
in an open
neighborhood of $\hat \lambda$, the value of
$\hat \iota$ is constant. We can further check that
$\hat \lambda_{\hat \iota} > \hat \lambda_{\hat\iota+1}$.
Therefore, at the associated $\hat x$,
we can employ \cite*[Theorem 3.1]{Tsing},
and we calculate
\[
\frac{\partial}{\partial x_j} \Gamma_s(F(\hat x))= \sum_{\ell=1}^{k} \frac{\partial  \phi_s(\lambda(F(\hat x)))}{\partial \lambda_{\ell}} h_j^{\ell}(\hat x)~,
\]
where
\[
h_j^{\ell}(\hat x)=\hat u_\ell^\top\frac{\partial F(\hat x)}{\partial x_j} \hat u_\ell = \hat u_\ell^\top F_{j\cdot}^\top F_{j\cdot} \hat u_\ell=(F_{j\cdot} \hat u_\ell)^2~.
\]

Calculating
\[
\frac{\partial \phi_s(\hat \lambda)}{\partial \lambda_\ell}  = \left\{\begin{array}{ll}
    1/\hat \lambda_\ell~, &\hbox{if $\ell\leq \hat \iota$;}\\[5pt]
     \frac{s-\hat \iota}{\sum_{i=\hat \iota+1}^{k}\hat \lambda_i}, &\hbox{if $\ell>\hat \iota$,}
    \end{array}\right.\\
\]
the result follows. \Halmos\endproof

Without the technical condition $\frac{1}{s-\hat \iota}\sum_{\ell=\hat \iota+1}^k \hat \lambda_{\ell}> \hat \lambda_{\hat\iota +1}$,
the formulae above still give a subgradient of $\Gamma_s$
(see \cite*{Weijun} for details).

Finally, considering  \ref{DDFact}, we can now present a short proof for Corollary \ref{cor:FactDoesntMatter}.
\medskip

\emph{Proof [Corollary \ref{cor:FactDoesntMatter}].}
As $	z_{{\tiny\mbox{DDFact}}}(C,s,A,b;F)=	z_{{\tiny\mbox{DFact}}}(C,s,A,b;F)$, it suffices to show that the objective value of \ref{DDFact} at any feasible solution $\hat x$, does not dependent on the  factorization $C=FF^\top$.
We have that   $F(\hat x):=F^\top\Diag(\hat x)F$ has the same non-zero eigenvalues as
the matrix
$\left(\Diag\left(\hat x\right)\right)^{\frac{1}{2}}FF^\top\left(\Diag\left(\hat x\right)\right)^{\frac{1}{2}}$, which is equal to $\left(\Diag\left(\hat x\right)\right)^{\frac{1}{2}}C\left(\Diag\left(\hat x\right)\right)^{\frac{1}{2}}$.
The result follows. \Halmos

\subsection{linx}

For $x\in[0,1]^n$ and $\gamma>0$, we define
$K_\gamma(x):=\gamma C\Diag(x)C+\Diag(\mathbf{e}-x)$.
The \emph{(scaled) linx  bound} for \ref{CMESP} is
	\[
	z_{{\tiny\mbox{linx}}}(C,s,A,b;\gamma):=\max \left\{
	\textstyle{\frac{1}{2}}(\ldet K_\gamma( x)-s\log \gamma)
	~:~ \mathbf{e}^{\top}x=s,~ Ax\leq b,~0\leq x\leq \mathbf{e} \right\}\tag{linx}\label{linx}
	\]
(see  \cite*{Kurt_linx}).

 It turns out that the linx bound is invariant under complementing (\cite*{Kurt_linx}; also see \cite*{FLbook}),
 while the factorization bound is not; therefore, we can obtain a different bound value
 for \ref{CMESP} by considering the factorization bound (but not the linx bound) on the complementary
 problem.

We can give an expression for the gradient and the Hessian of the \ref{linx} objective function.
Using well-known facts, we can work out that
for all $\hat{x}$ in the domain of the objective function, we have
\[\nabla f(x) = \frac{1}{2}\left(\diag(\gamma CK_\gamma(\hat x)^{-1} C)-\diag(K_\gamma(\hat x)^{-1})\right),
\]
and
\begin{align*}
&\nabla^2\left(\textstyle{\frac{1}{2}}\ldet K_\gamma( \hat x)\right) =\\
&\quad\textstyle{\frac{1}{2}}\Bigg(- \gamma^2 (C K_\gamma(\hat x)^{-1} C)\circ (CK_\gamma(\hat x)^{-1} C)  \\
&\quad\quad + \gamma \Big((K_\gamma(\hat x)^{-1} C)\circ (K_\gamma(\hat x)^{-1} C)
 +  ((K_\gamma(\hat x)^{-1} C)\circ (K_\gamma(\hat x)^{-1} C))^\top\Big)
 - K_\gamma(\hat x)^{-1}\circ K_\gamma(\hat x)^{-1}\Bigg).\\
\end{align*}


 \subsection{Mixing}\label{sec:mixing}

We consider  $m\geq 1$ convex relaxations for \ref{CMESP}, indexed by $i=1,\ldots,m$:
\[
v_i:=
\max\left\{ f_i(L_i(x)) ~:~
\mathbf{e}^\top x = s,~ Ax\leq b,~ 0\leq x\leq \mathbf{e}\right\},
\]
where, for $i=1,\ldots,m$,  $k_i\leq n$,   $L_i:\mathbb{R}^n\rightarrow\mathbb{S}^{k_i}_+$
are affine functions,  and $ f_i:\mathbb{S}^{k_i}_+\rightarrow \mathbb{R}$
are concave functions.
We write  $L_i(x):=L_{i0}+L_{i1}x_1+\cdots+L_{in}x_n$ and $L_{ij}\in\mathbb{S}^{k_i}$,   for $i=1,\ldots,m$ and $j=0,\ldots,n$, and we
note that the objective functions of \ref{DDFact}, complementary \ref{DDFact}, and \ref{linx} can  be written as $f_i(L_i(x))$ (see \S\ref{sec:details}).


 For a ``weight vector'' $\alpha\in\mathbb{R}^m_+$, such that  $\mathbf{e}^\top \alpha=1$,   we
define the \emph{mixing bound} (see \cite*{Mixing} for a more general setting):
\begin{equation}\label{prob:mixing}\tag{mix}
v(\alpha) := \max\left\{ \sum_{i=1}^m \alpha_i f_i(L_i(x)) ~:~
\mathbf{e}^\top x = s,~ Ax\leq b,~ 0\leq x\leq \mathbf{e}
\right\}.
\end{equation}
The goal is to minimize the mixing bound over  $\alpha$ (and any parameters for the individual bounds).


We construct the Lagrangian dual of \ref{prob:mixing} for a broad class of cases that
covers our applications of mixing.
For $i=1,\ldots,m$,   we assume  that for any given $\hat\Theta_i\in\mathbb{S}^{k_i}_{++}$, there is a closed-form solution  $\hat{W_i}$ to  $\sup\{  f_i(W_i)-\hat\Theta_i\bullet W_i~:~W_i\succeq 0\}$, such that $\hat\Theta_i\bullet \hat{W}_i=:\rho_i\in\mathbb{R}$ and   $\Omega_i:\mathbb{S}^{k_i}_{++}\rightarrow\mathbb{R}$, is defined by $\Omega_i(\hat\Theta_i)
:= f_i(\hat{W}_i)$. Furthermore, we assume that the supremum is $+\infty$ if $\hat \Theta_i\nsucc 0$.


With the assumptions above, the Lagrangian dual problem of \ref{prob:mixing} is equivalent to (see the Appendix for details)

\[
\begin{array}{ll}
&	z_{{\tiny\mbox{Dmix}}}(C,s,A,b):=\\
&~~ 	\min~  \displaystyle\sum_{i=1}^m \alpha_i\Big( \Omega_i(\Theta_i) - \rho_i+ \Theta_i\bullet L_{i0}\Big)
+ \nu^\top \mathbf{e}  + \pi^\top b +\tau s\\
     &~~ \mbox{subject to:}\\
     &\quad   \displaystyle\sum_{i=1}^m \alpha_i\Big(\Theta_i\bullet L_{ij }\Big) + \upsilon_j  - \nu_j   - \pi^\top A_{\cdot j} - \tau=0,~\mbox{ for } j \in N,\\
&\quad \Theta_1\succ 0,\ldots,\Theta_m\succ 0, ~\upsilon\geq 0, ~\nu\geq 0, ~\pi\geq 0.
\end{array}
\tag{Dmix}\label{Dmix}
\]

\begin{thm} \label{thm:fixmix}
 Let
\begin{itemize}
    \item $\mbox{LB}$ be the objective-function value of a feasible solution for \ref{CMESP},
    \item $(\hat{\Theta}_1,\ldots,\hat{\Theta}_m,\hat{\upsilon},\hat{\nu},\hat{\pi},\hat{\tau})$ be a feasible solution for \ref{Dmix} with objective-function \hbox{value $\hat{\zeta}$.}
\end{itemize}
Then, for every optimal solution  $x^*$ for \ref{CMESP}, we have:
\[
\begin{array}{ll}
x_j^*=0, ~ \forall ~ j\in N \mbox{ such that } \hat{\zeta}-\mbox{LB} < \hat{\upsilon}_j~,\\
x_j^*=1, ~ \forall ~ j\in N \mbox{ such that } \hat{\zeta}-\mbox{LB} < \hat{\nu}_j~.\\
\end{array}
\]
\end{thm}
\proof{Proof.}
Analogous to the proof of Theorem \ref{thm:fixFact}.
\Halmos\endproof

Next, we generalize to \ref{Dmix}, the procedure presented in  \S\ref{sec:DDFact} to construct a feasible solution of \ref{DFact} from a feasible solution  of \ref{DDFact}.  A good feasible solution for \ref{Dmix} can be used to validate the quality of the solution obtained for  \ref{prob:mixing}, and to fix variables by applying the result of Theorem \ref{thm:fixmix}.

We let $\hat x$ be  a  feasible solution  of \ref{prob:mixing}  in the domain of $f_i$ and  define   $\hat{W}_i:=L_i(\hat x)$, for  $i=1,\ldots,m$. First, we assume that  it is possible to compute $\hat \Theta_i$, such that $\Omega_i(\hat{\Theta}_i)=f_i(\hat W_i)$,  for $i=1,\ldots,m$.
Then, the minimum duality gap between $\hat x$ in \ref{prob:mixing} and feasible solutions
of \ref{Dmix} of the form $(\hat{\Theta}_1,\ldots,\hat{\Theta}_m,{\upsilon},{\nu},{\pi},{\tau})$
is the optimal value of the linear program

\begin{equation*}\label{mingapprobmix}\tag{$G(\hat{\Theta}_1,\ldots,\hat{\Theta}_m)$}
\begin{array}{l}
\min~
 \nu^\top \mathbf{e}  + \pi^\top b +\tau s - \displaystyle \sum_{i=1}^m \alpha_i\Big( \rho_i- \hat{\Theta}_i\bullet L_{i0}\Big)\\
      \quad \mbox{subject to:}\\
     \quad   \upsilon_j  - \nu_j   - \pi^\top A_{\cdot j} - \tau =- \displaystyle \sum_{i=1}^m \alpha_i\left(\hat{\Theta}_i\bullet L_{ij }\right),~\mbox{ for } j \in N,\\
\quad \upsilon\geq 0, ~\nu\geq 0, ~\pi\geq 0.
\end{array}
\end{equation*}

\noindent
Analogously to  \ref{mingapprob}, we can verify that \ref{mingapprobmix} is always feasible and has a simple closed-form solution
for \ref{MESP}. In fact, the only differences between
\ref{mingapprob} and \ref{mingapprobmix} are the constant in the objective function
and the right-hand side of the constraints.

\subsubsection{Considering  \ref{DDFact}, complementary \ref{DDFact}, and \ref{linx} in \ref{prob:mixing} }\label{sec:details}

Considering $f_i(L_i(x))$  as the objective function of \ref{DDFact} we have $f_i(\cdot):=\Gamma_s(\cdot)$ and $L_i(x):=F^\top \Diag(x)F$, so $k_i:=k$, $L_{i0}:=0$ and $L_{ij}:=F_{j\cdot}^\top F_{j\cdot}$ , for $j=1,\ldots,n$.
We also  have  $\rho_i:=s$ and $\Omega_i(\Theta_i):=-\sum_{\ell=k-s+1}^k \log(\lambda_\ell(\Theta_i))$. For a given   feasible solution  $\hat x$  of \ref{prob:mixing}  in the domain of $f_i$ and    $\hat{W}_i:=L_i(\hat x)$, we construct  $\hat \Theta_i$ as discussed in \S\ref{sec:DDFact}, and we see in \eqref{sumbeta}, that  $\Omega_i(\hat{\Theta}_i)=f_i(\hat W_i)$.

Considering $f_i(L_i(x))$  as the objective function of complementary \ref{DDFact} we have $f_i(\cdot):=\Gamma_{n-s}(\cdot) + \ldet C$ and
$L_i(x):=F^{-1} \Diag(\mathbf{e}-x)F^{-\top}$, so $k_i:=k(=n)$, $L_{i0}:=F^{-1} F^{-\top}$ and $L_{ij}:=-F_{j\cdot}^{-1} F_{j\cdot}^{-\top}$ , for $j=1,\ldots,n$.
We also have  $\rho_i:=n-s$ and $\Omega_i(\Theta_i):=-\sum_{\ell=k-n+s+1}^k \log(\lambda_\ell(\Theta_i))+\ldet C$. For a given   feasible solution  $\hat x$  of \ref{prob:mixing}  in the domain of $f_i$ and    $\hat{W}_i:=L_i(\hat x)$, we construct  $\hat \Theta_i$ as discussed in \S\ref{sec:DDFact}, and we see in \eqref{sumbeta}, that  $\Omega_i(\hat{\Theta}_i)=f_i(\hat W_i)$.

 Considering $f_i(L_i(x))$  as the objective function of \ref{linx} we have $f_i(\cdot):=\frac{1}{2}\left(\ldet(\cdot) -s\log\gamma\right)$ and $L_i(x):=\gamma C\Diag(x)C+\Diag(\mathbf{e}-x)$, so  $k_i:=n$, $L_{i0}:=I$ and $L_{ij}:=\gamma C_{j\cdot}^\top C_{j\cdot}- \mathbf{e}_{j}\mathbf{e}_{j}^\top$ , for $j=1,\ldots,n$. We also have $\rho_i:=n/2$ and $\Omega_i(\Theta_i):=-\frac{1}{2}\left(\ldet(2\Theta_i)+s\log\gamma\right)$. For a given   feasible solution  $\hat x$  of \ref{prob:mixing}  in the domain of $f_i$ and    $\hat{W}_i:=L_i(\hat x)$, we set  $\hat \Theta_i:=\frac{1}{2}\hat{W}_i^{-1}$, and we see that  $\Omega_i(\hat{\Theta}_i)=f_i(\hat W_i)$.

We note that if the objective of \ref{prob:mixing} is a weighted combination of the three functions mentioned above and  $\hat x$ is an optimal solution to \ref{prob:mixing}, then the optimal objective value of \ref{mingapprobmix} is zero, that is, the dual solution constructed to \ref{Dmix} is also optimal.

\section{Implementation and experiments}\label{sec:impl}

\subsection{Setup for the computational  experiments}
\cite*{Weijun} worked with solving \ref{DDFact} with respect to \ref{MESP}, using
a custom-built Frank-Wolf (see \cite*{FrankWolfe}) style
code, written in Python. They only worked with the relaxation,
and did not seek to solve \ref{MESP} to optimality.

The linx bound for \ref{CMESP} was introduced by \cite*{Kurt_linx}, where bound calculations
were carried out with the conical-optimization software \texttt{SDPT3} (see \cite*{SDPT3}), within
the very-convenient \texttt{Yalmip} \texttt{Matlab} framework (see \cite*{Yalmip}), and a full branch-and-bound code for \ref{MESP} was written in \texttt{Matlab}.

In our experiments, we calculate all of our bounds using a single state-of-the-art commercial nonlinear-programming solver, to
facilitate fair comparisons between bounding methods, and also to see what
is possible in such a computational setting. 



We experimented on instances of \ref{MESP}  and  \ref{CMESP} with \ref{linx},   \ref{DDFact} and complementary \ref{DDFact} (i.e,
\ref{DDFact} applied to \ref{CMESP-comp}).
We ran our experiments under Windows, on an Intel Xeon E5-2667 v4 @ 3.20 GHz processor equipped with 8 physical cores (16 virtual cores) and 128 GB of RAM.
We implemented our code in \texttt{Matlab} using the commercial  software \texttt{Knitro}, version 12.4, as our nonlinear-programming solver.
\texttt{Knitro} offers BFGS-based algorithms and  Newton-based algorithms to solve nonlinear programs. In the first case, \texttt{Knitro} only needs function
values and gradients from the user, in the latter, \texttt{Knitro} also needs second derivatives. By experimenting on top of one state-of-the-art
general-purpose nonlinear-programming code,
we hoped to get good and rapid convergence
and get running times that can reasonably be compared for the different relaxations.
In all of our experiments we set \texttt{Knitro} parameters\footnote{see \url{https://www.artelys.com/docs/knitro/2_userGuide.html}, for details} as follows:
\verb;algorithm;~$=3$ to use an active-set method,
\verb;convex;~$=1$ (true),
\verb;gradopt;~$=1$ (we provided exact gradients),
\verb;maxit;~$=1000$.  We set
\verb;opttol;~$=10^{-10}$, aiming to satisfy the  KKT optimality conditions to a very tight tolerance. We set
 \verb;xtol;~$=10^{-15}$ (relative tolerance for lack of progress in the solution point) and
 \verb;feastol;~$=10^{-10}$ (relative tolerance for the feasibility error), aiming for
 the best solutions that we could reasonably find.


\subsection{Test instances}

To compare the bounds obtained with the three relaxations, we consider four covariance matrices  from the literature, with  $n=63, 90, 124, 2000$.  For each matrix, we consider different values of $s$ defining a set of test instances of \ref{MESP}. The  $n=63$ and $n=124$ matrices are benchmark covariance matrices
obtained from J. Zidek (University of British Columbia), coming from an application to re-designing an environmental monitoring network;
see \cite*{Guttorp-Le-Sampson-Zidek1993} and \cite*{HLW}.
The $n=90$ matrix is based on temperature data from monitoring stations
in the Pacific Northwest of the United States; see \cite*{Kurt_linx}.
These $n=63, 90, 124$ matrices are all nonsingular.
All of these matrices have been used extensively in testing and developing algorithms for MESP; see \cite*{KLQ,LeeConstrained,AFLW_Using,LeeWilliamsILP,HLW,AnstreicherLee_Masked,BurerLee,Anstreicher_BQP_entropy,Kurt_linx}. The largest covariance matrix that we considered in our experiments is an $n=2000$ matrix with rank $949$, based on Reddit data, used in \cite*{Weijun} and from  \cite*{Dey2018} (also see \cite*{Munmun}).


To ameliorate some instability in running times,
 for $n=63, 90, 124$ we repeated every experiment ten times, and for
 $n=2000$, we  repeated every experiment five times, and present average
 timing results.

 \subsection{Numerical experiments for $n=63,90,124$}

For the three nonsingular covariance matrices used in our experiments, we solved  \ref{linx},   \ref{DDFact} and complementary \ref{DDFact},  for all $2\leq s\leq n-1$. For each matrix, we present four plots. In the first plots of Figures \ref{fig:data63}, \ref{fig:data90} and \ref{fig:data124}, we present the integrality gap for each bound and each $s$.  Each such gap is
given by the difference between the upper bound computed by solving the relaxation  and a lower bound obtained using a heuristic of \cite*[Sec. 4]{LeeConstrained} followed by a simple local search (see \cite*[Sec. 4]{KLQ}).

\begin{figure}[htbp]
	\centering
	\scalebox{0.65}{
		\includegraphics{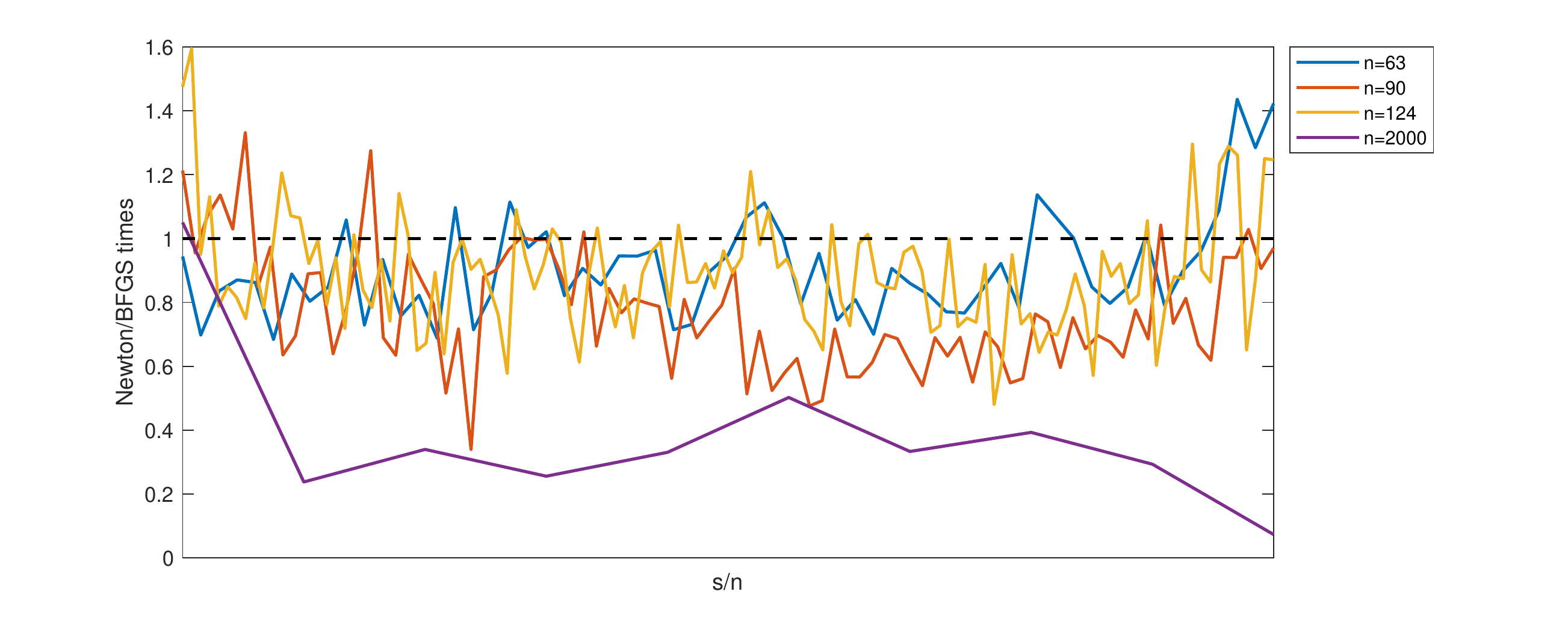}
	}
	\caption{Newton/BFGS time for \ref{linx}}
	\label{fig:Linx:time_ratio}
\end{figure}

In the second plots of  those figures, we present the average wall-clock times (in seconds) used by \texttt{Knitro} to solve the relaxations.  Some observations about the  times presented are important. First, we note that the times depicted on the plots correspond to the application of a BFGS-based algorithm to solve \ref{DDFact} and complementary \ref{DDFact}, and to the application of a Newton-based algorithm to solve \ref{linx}. As an experiment,  we also applied a BFGS-based algorithm to solve \ref{linx}, not passing the Hessian to the solver, but, as expected, the results were worse concerning both time and convergence of the algorithm. The difference between the times can be seen in Table \ref{tab:wallclocktime} (aggregated over $s$) and Figure \ref{fig:Linx:time_ratio}.  On the other hand,  we did not apply a Newton-based algorithm to \ref{DDFact} and complementary \ref{DDFact} because we cannot guarantee that the objective function of these relaxations
is differentiable at every iterate  (of the nonlinear-programming solver).  We should also note that the times shown in the plots for \ref{linx} do not include the times to compute the value of the parameter $\gamma$ in the problem formulation. This parameter value has a great impact on
the \ref{linx} bound. We present the times to compute them, aggregated over $s$, in Table \ref{tab:wallclocktime}.
Finally, we should mention that \texttt{Knitro} did not prove optimality for several instances solved.
However, we could confirm the optimality of \emph{all} solutions returned by \texttt{Knitro}, up to the optimality tolerance considered, by constructing a dual solution with duality gap less than the tolerance, with respect to the primal solution obtained by \texttt{Knitro}. To construct the dual solutions, we solved various special cases of the linear program  \ref{mingapprobmix} (see \S\ref{sec:mixing}). For \ref{DDFact},
  the construction    uses \eqref{defbeta},
where we took $\epsilon=0$, which  gives us
a dual solution that is feasible within numerical accuracy.



  On our experiments with instances of \ref{MESP} (i.e., no side constraints), these linear programs have closed-form solutions and the times to compute them are not significant.

In  the third plots of Figures \ref{fig:data63},  \ref{fig:data90} and \ref{fig:data124}, we demonstrate the capacity of the mixing methodology described in  \S\ref{sec:mixing} to decrease the integrality gap.  As observed in \cite*{Mixing}, the methodology is particularly effective when considering in \ref{prob:mixing}, a weighted sum of the objective functions of two relaxations, such that the bounds obtained by each relaxation are close to each other. We exploit this observation in our experiments. For each covariance matrix, we select one or more pairs of relaxations for which the integrality-gap curves (presented in the first plots of the figures) cross each other at some point. Then, we mix these two relaxations and compute new mixed bounds for all values of $s$ in a promising interval, approximately centered at the point where the two curves cross. To select the parameter $\alpha$ that weights the objective in \ref{prob:mixing}, we  simply apply a bisection algorithm.

Finally, in the fourth plots of Figures \ref{fig:data63}, \ref{fig:data90}, and \ref{fig:data124}, we demonstrate  how effective the strategy described in Theorems \ref{thm:fixFact} and \ref{thm:fixmix} can be to fix variables (e.g., the context could be fixing variables at the root node of the enumeration tree in applying a branch-and-bound algorithm).  In all of our experiments, we use a fixing threshold of $10^{-10}$ which can be considered as
rather safe in the context of the accuracy that
we use to compute the relevant quantities.
Although the mixing strategy can decrease the integrality gap for some instances, in our experiments this improvement is not enough to allow more variables to be fixed. Therefore, we do not consider the mixed bounds in these plots.

\subsection{Analysis of the results for $n=63,90,124$}

The analysis of the plots for the $n=63$ and $n=90$  covariance matrices  are very similar to each other and are summarized in the following.

\begin{figure}[htbp]
	\centering

	\scalebox{0.55}{
		\includegraphics{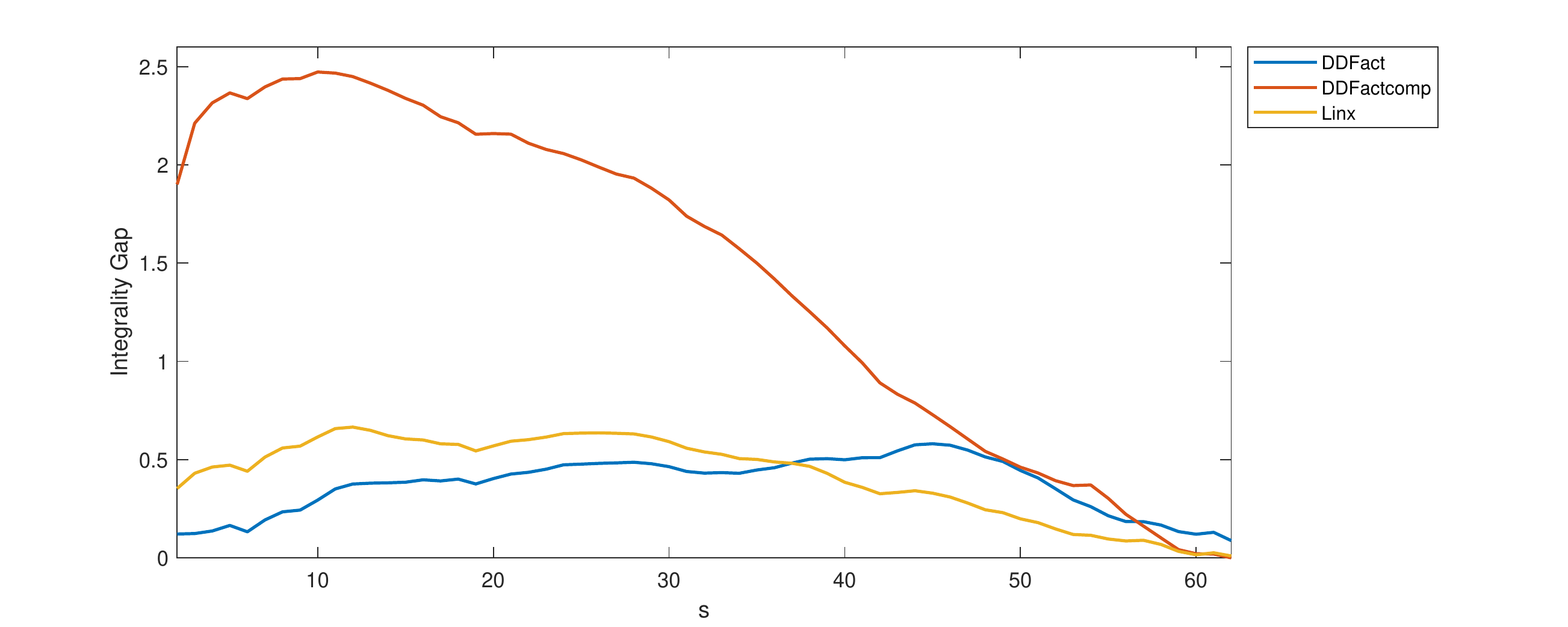}}
	\scalebox{0.55}{
	\includegraphics{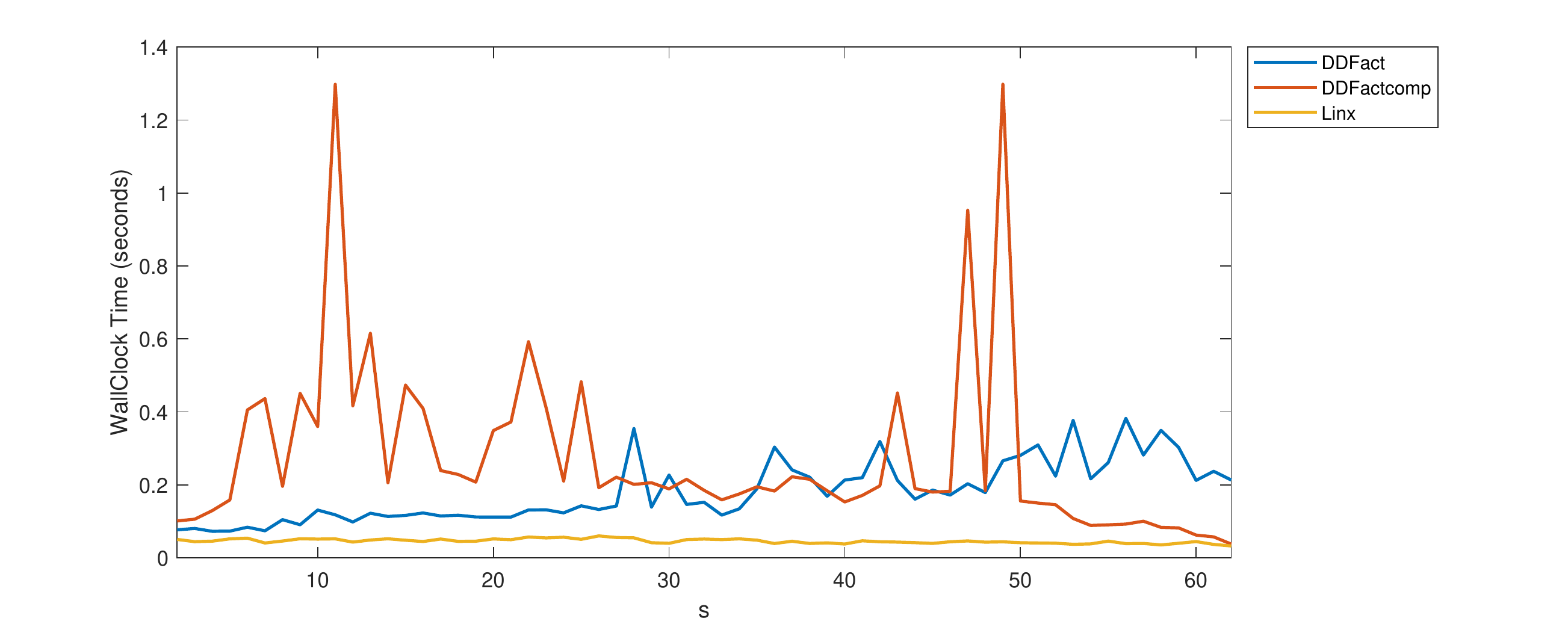}
	}
	\scalebox{0.55}{
	\includegraphics{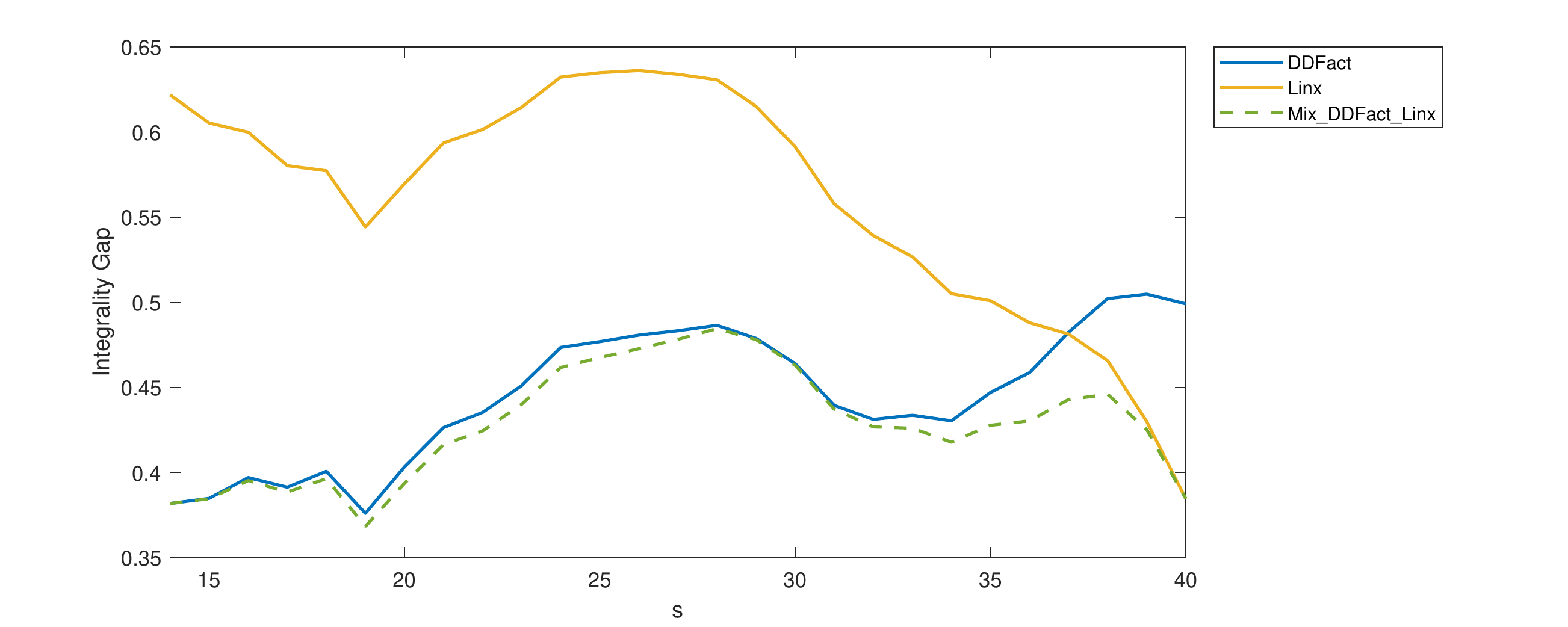}}
	\scalebox{0.55}{
	\includegraphics{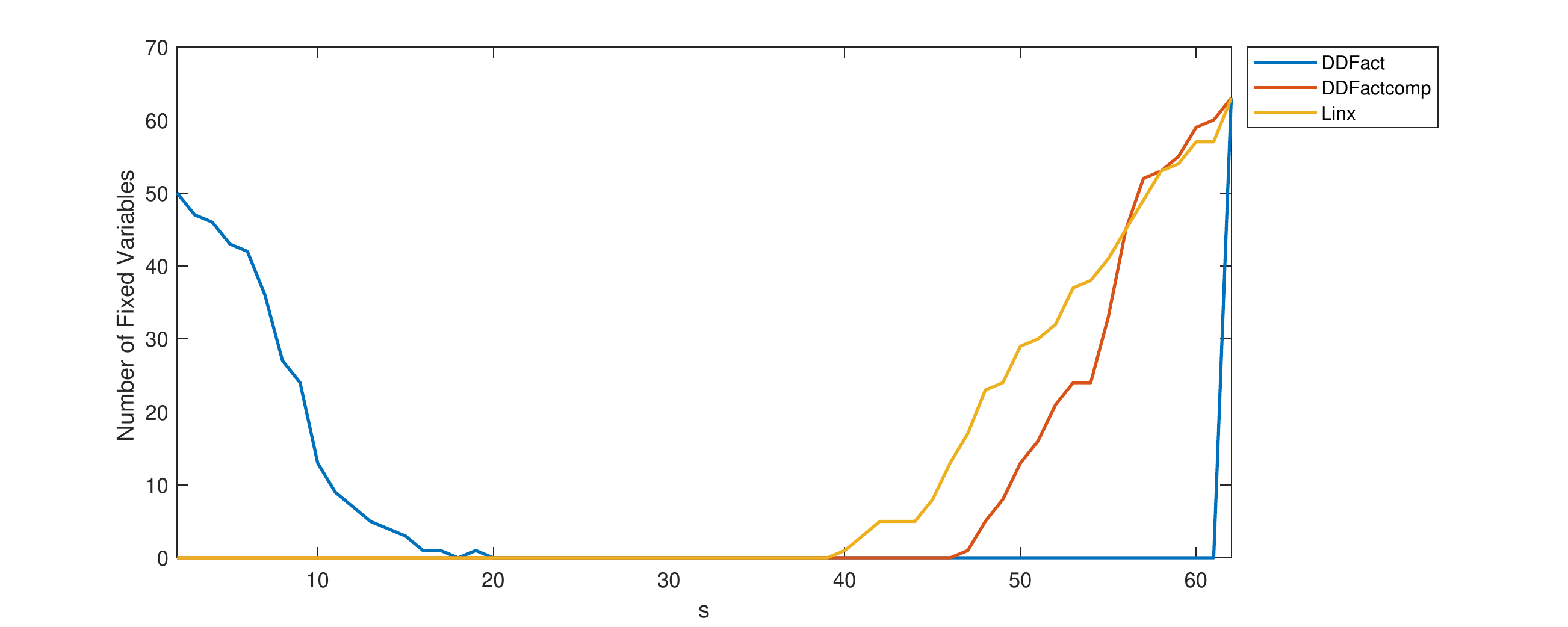}
	}
	\caption{Bounds/times comparison and effect of the mixing  and  variable-fixing methodologies for $n=63$}
	\label{fig:data63}
\end{figure}

\begin{figure}[htbp]
	\centering

	\scalebox{0.55}{
		\includegraphics{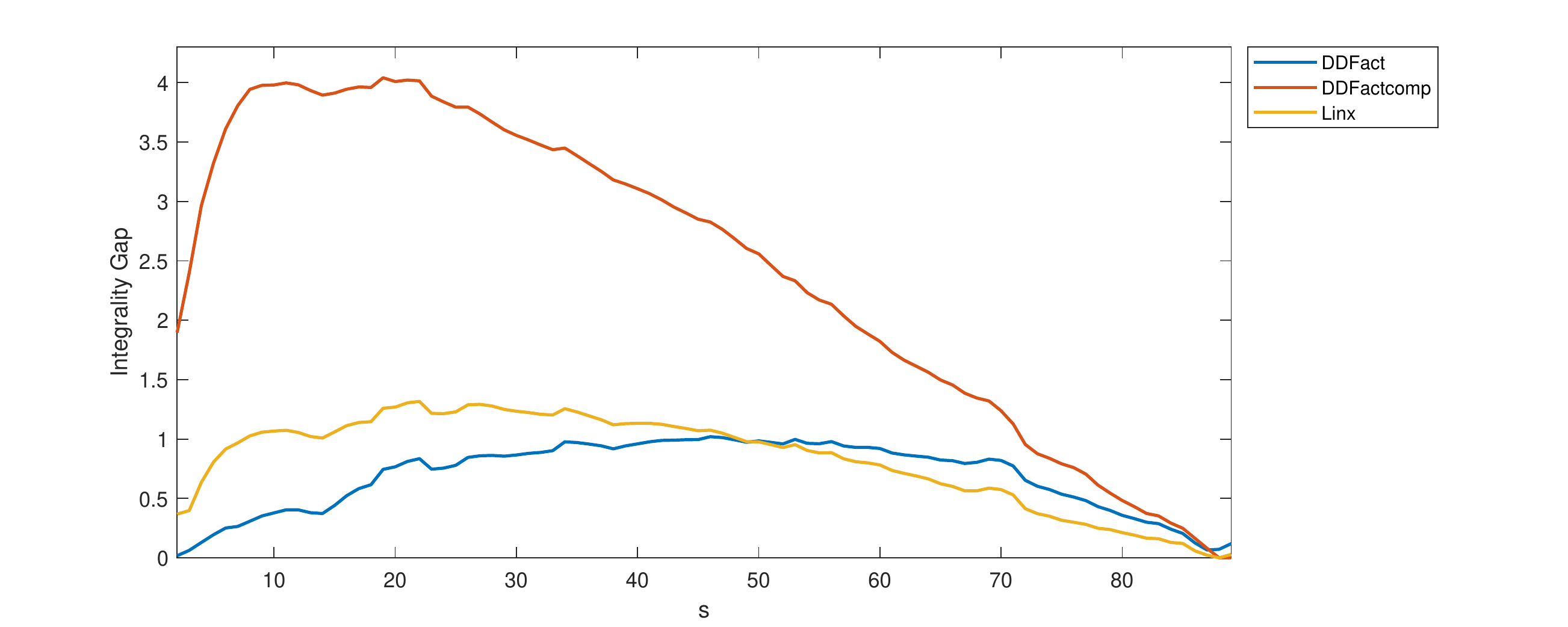}}
	\scalebox{0.55}{
	\includegraphics{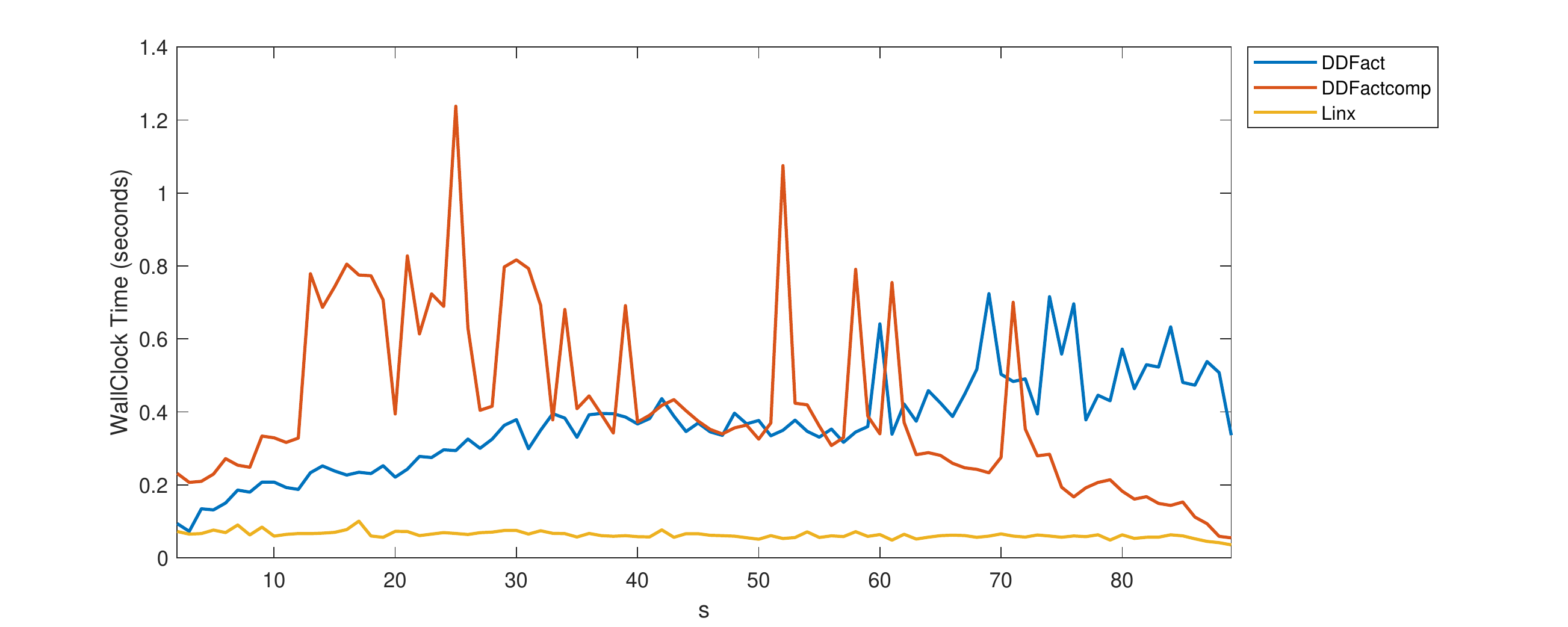}
	}
	\scalebox{0.55}{
	\includegraphics{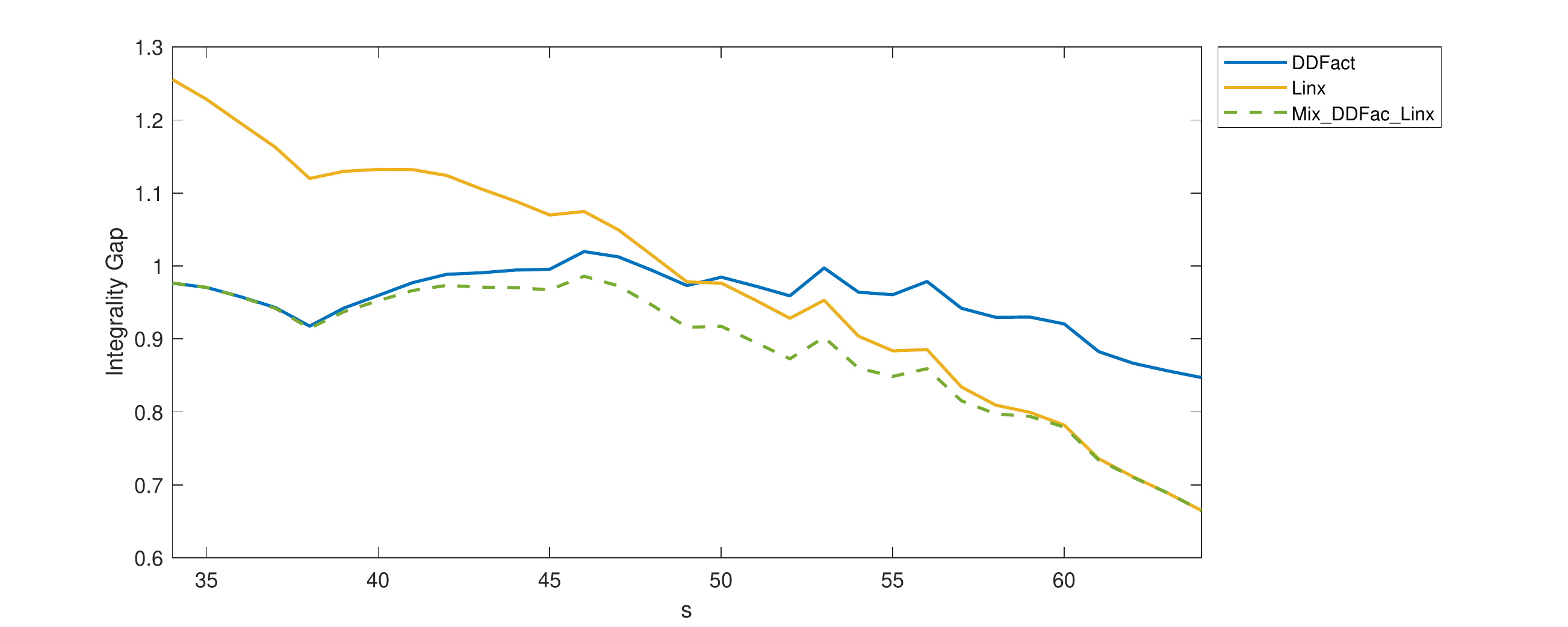}}
	\scalebox{0.55}{
	\includegraphics{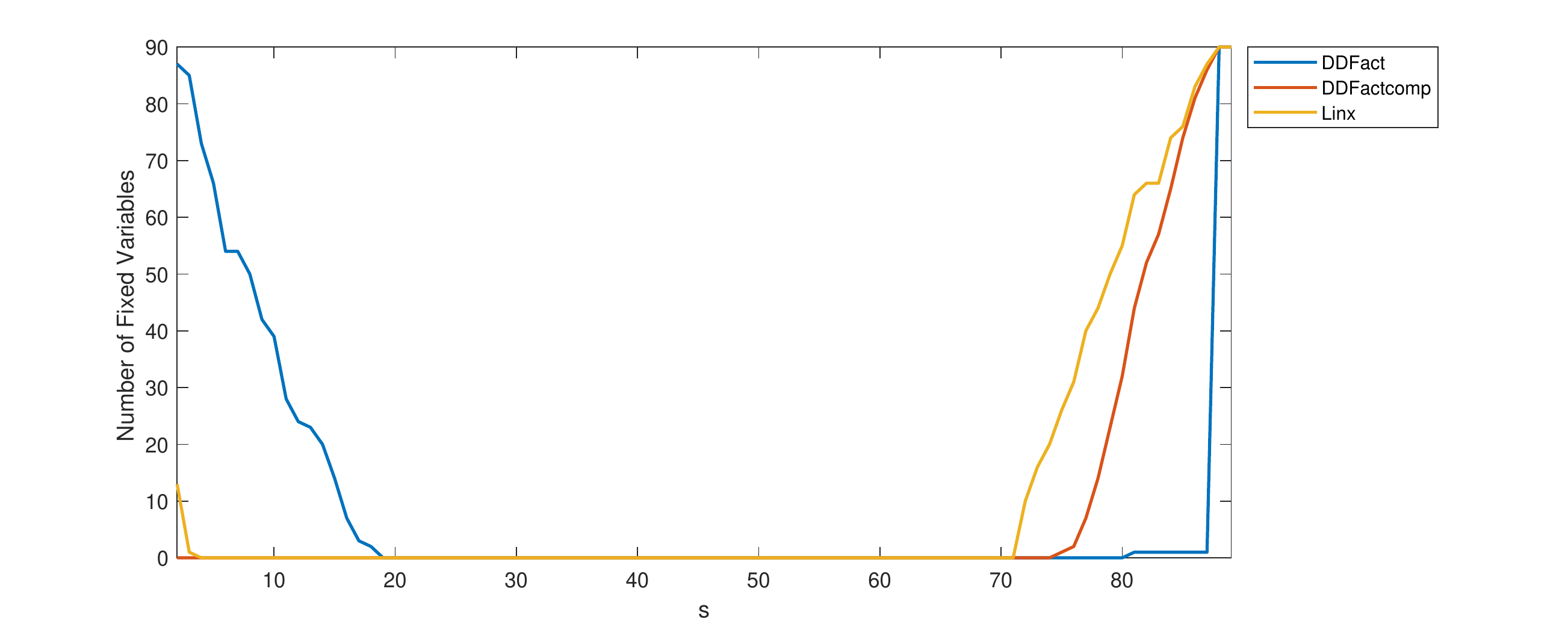}
	}
	\caption{Bounds/times comparison and effect of the mixing  and  variable-fixing methodologies for $n=90$}
	\label{fig:data90}
\end{figure}

\begin{itemize}
\item  We see from the first plots of Figures \ref{fig:data63} and \ref{fig:data90} that the complementary \ref{DDFact} bound is not competitive with the \ref{DDFact} and \ref{linx} bounds for these instances. For most values of $s$ the first bound is much worse than the two others. The complementary \ref{DDFact} bound is only a bit better than the \ref{DDFact} bound for very large values of $s$ and it is never better than the \ref{linx} bound.  The integrality-gap curves for \ref{DDFact} and \ref{linx}  cross at points close to an intermediate value of  $s$. For smaller  $s$, \ref{DDFact} gives the  best bound and for larger $s$, \ref{linx} is the winner.

\item Concerning the wall-clock time to compute the bounds, we see again a big disadvantage of complementary \ref{DDFact} in the second plots of Figures \ref{fig:data63} and \ref{fig:data90}.  Although it is faster than \ref{DDFact} on some instances with large $s$, we see that on most instances its time is much longer than the times for the  two other relaxations, and with the greatest variability among  the different values of $s$. The solution of  \ref{linx} is always significantly faster than the solution of the  two other relaxations for these instances. Finally, we  note that the variation in the time to solve \ref{linx} for all values of $s$ is less than $1\%$, while there is a great variability for the other two relaxations.
\item   We see in the first plots of Figures \ref{fig:data63} and \ref{fig:data90} that the curves corresponding to \ref{DDFact} and \ref{linx} cross at points close to  intermediate values of $s$, indicating a promising interval of values to mix these relaxations for both  $n=63$ and $n=90$.   Considering such intervals, we see in the third plots of Figures \ref{fig:data63} and \ref{fig:data90}, how mixing \ref{DDFact} and \ref{linx} can in fact, decrease the integrality gap for some instances, being mostly effective for the values of $s$ for which the \ref{DDFact} bound and the \ref{linx} bound are very close to each other.
\item In the fourth plots of Figures \ref{fig:data63} and \ref{fig:data90} we verify the increasing capacity to fix variables as the bound gets stronger.  Interestingly, we see that for large values of $s$,  complementary \ref{DDFact} bounds can lead to more variables fixed than the better \ref{DDFact} bound.
\end{itemize}

\begin{figure}[htbp]
	\centering

	\scalebox{0.55}{
		\includegraphics{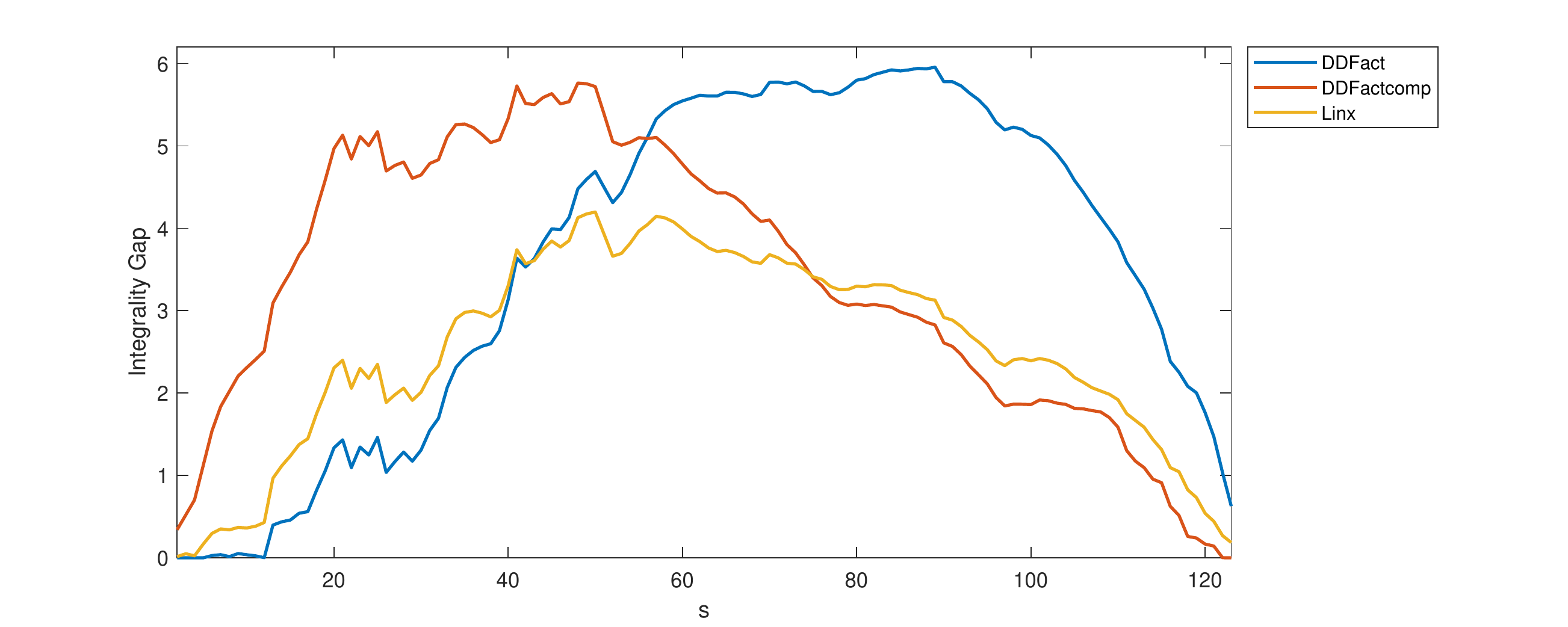}}
	\scalebox{0.55}{
	\includegraphics{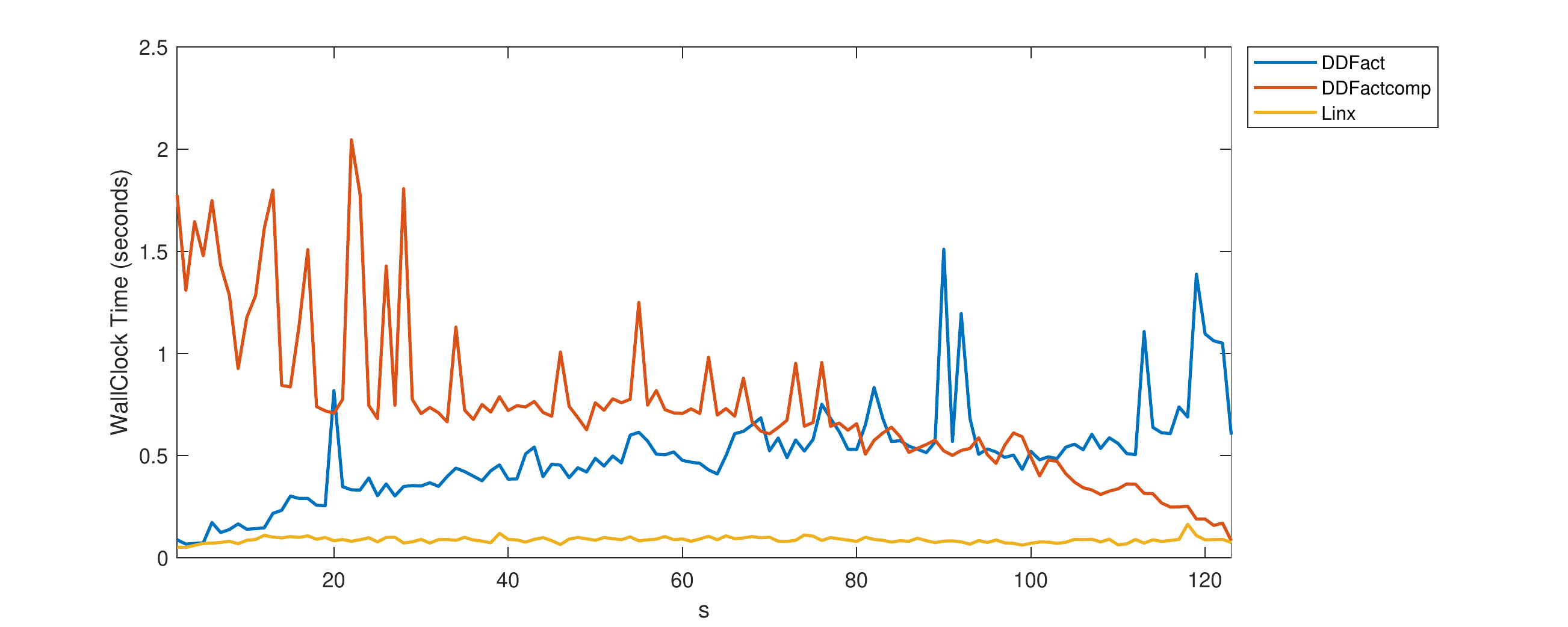}
	}
	\scalebox{0.60}{\phantom{xxx}
	\includegraphics{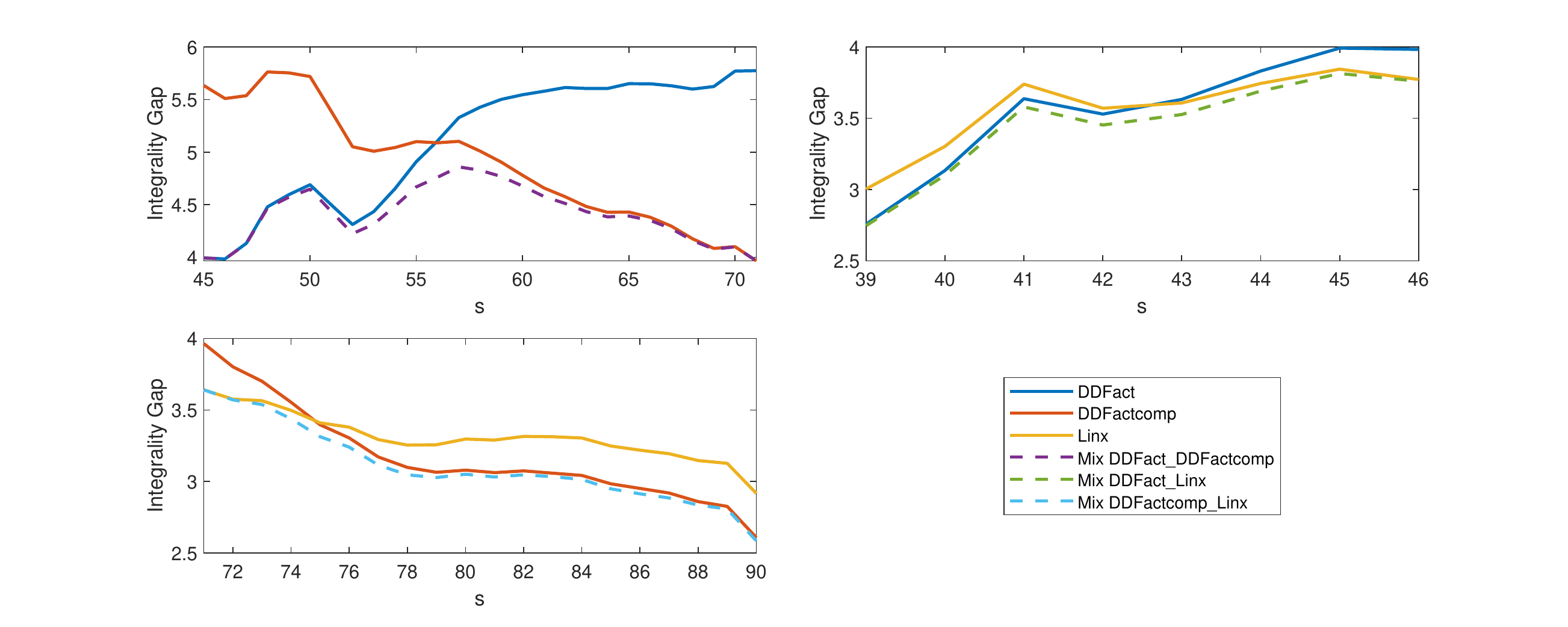}}
	\scalebox{0.55}{
	\includegraphics{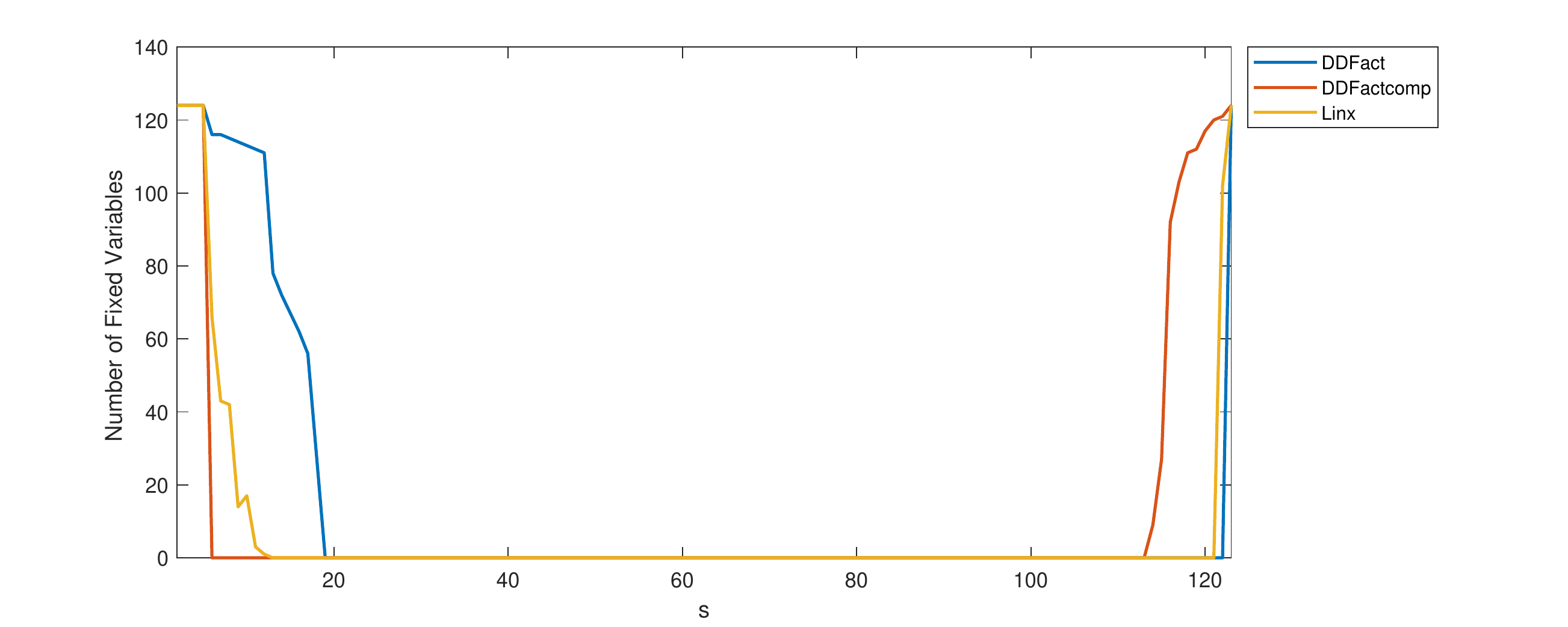}
	}
	\caption{Bounds/times comparison and effect of the mixing  and  variable-fixing methodologies for $n=124$}
	\label{fig:data124}
\end{figure}

For $n=124$, we have a slightly different analysis because, as we see in the first plot of Figure \ref{fig:data124}, the complementary \ref{DDFact} bound becomes better than the \ref{DDFact} bound for all  $s$ larger than  an  intermediate value.  We note that we can observe this  same behavior with the ``NLP" relaxation for \ref{CMESP} used in \cite*{AFLW_Using}.  Moreover, we see in the first plot of Figure \ref{fig:data124}, two points where the curves cross, showing three interesting intervals for $s$, where each one of the three relaxations gives the best bound. Concerning the wall-clock time, the observations about the second plot of  Figure \ref{fig:data124} are similar to the ones about  the second plots of Figures \ref{fig:data63} and \ref{fig:data90}, confirming that the time to solve  \ref{linx} is shorter  and with a smaller variability with $s$, when compared to the two other relaxations. In the third plot of Figure \ref{fig:data124},  we exploit the three crossing points of the integrality-gap curves for $n=124$, and show separately the capacity of the mixing methodology to decrease the gaps when we mix the two relaxations corresponding to each crossing point. 
It is interesting to note that the mixing methodology is more effective when the  crossing curves are less flat at those points, that is, when the gaps change faster as  $s$ changes.   Finally, we have  different observations about the fourth plot of Figure \ref{fig:data124} when compared to the smaller instances, concerning the capacity of the relaxations to fix variables.  As \ref{DDFact} and complementary \ref{DDFact} lead to very small integrality gaps at both ends of the curves, we observe in the fourth plot of Figure \ref{fig:data124} their stronger capacity of fixing variables on the corresponding values of $s$, when compared to  \ref{linx}. For $n=124$, we see that \ref{linx} gives the best bounds for intermediate values of $s$ only. These are  clearly the most difficult instances for $n=124$. Therefore, the integrality gap is usually not small enough on these instances to allow variable fixing.

\subsection{Numerical experiments with the large instance ($n=2000$)}

\begin{figure}[htbp]
	\centering

	\scalebox{0.55}{
		\includegraphics{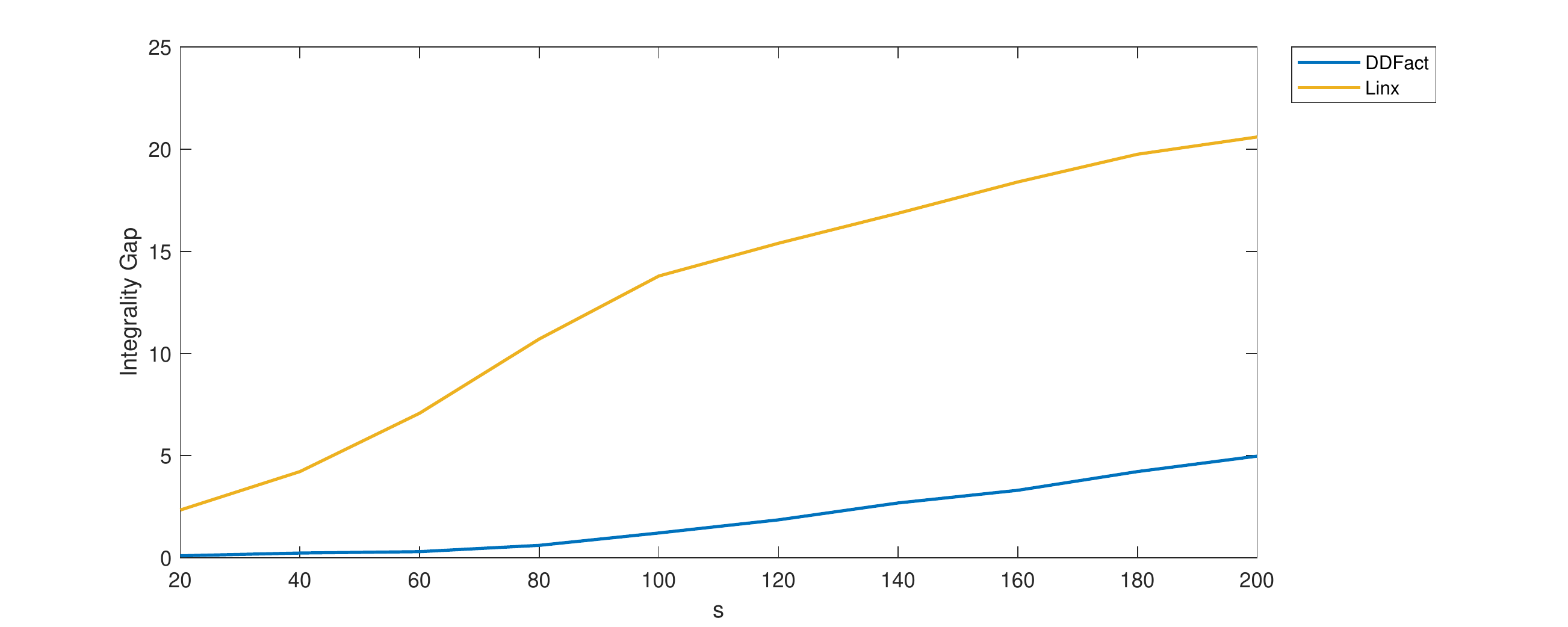}}
	\scalebox{0.55}{
	\includegraphics{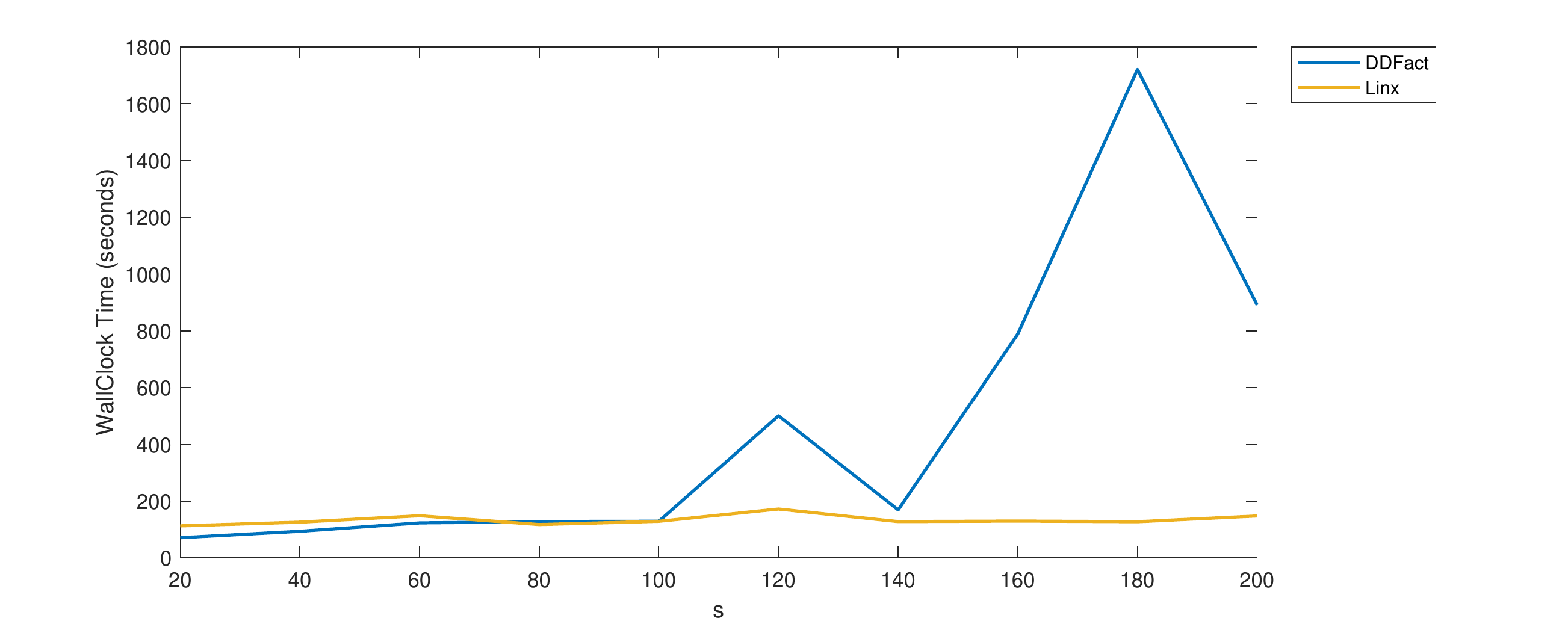}
	}
	\scalebox{0.55}{
	\includegraphics{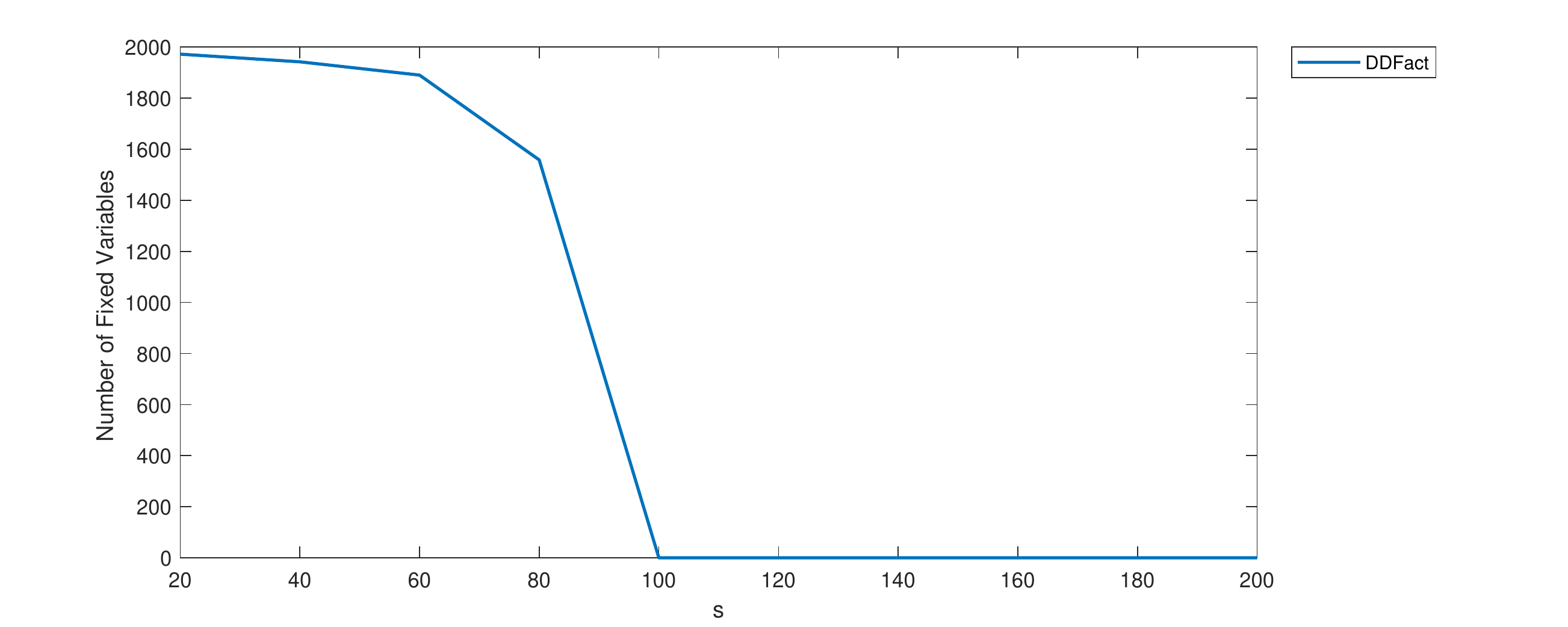}
	}
	\caption{Bounds/times comparison and effect of the  variable-fixing methodology for $n=2000$}
	\label{fig:data2000}
\end{figure}

The bounds computed for the $n=2000$ matrix are analyzed in Figure \ref{fig:data2000}. As this larger matrix is singular, we could not apply the complementary \ref{DDFact} relaxation to obtain a bound.
In the first plot, we present the integrality gaps for \ref{DDFact} and  \ref{linx}  for all $20\leq s\leq 200$ that are multiples of 20.
For lower bounds in computing the gaps, we obtained them by the same heuristic applied to our smaller instances with $n=63,90,124$.
We clearly see the superiority of \ref{DDFact} for this input matrix, for these relatively small values of $s$, following the behavior observed for the smaller instances. Concerning the wall-clock time, we still see in the second plot  that \ref{linx} can be solved faster on the most difficult instances with $s\geq 100$, and once more, we see a very small variability in the times for \ref{linx}, unlike what we see for \ref{DDFact}. The significant rank deficiency of the covariance matrix of these instances would seem to be a disadvantage for \ref{linx} relative to \ref{DDFact}, with regard to the computational time needed to solve them;
this is because the order of the matrix considered in the objective function of \ref{linx} is always equal to the order of the covariance matrix, while for \ref{DDFact} it is given by its rank. As \ref{linx} could not fix variables for any value of $s$, we present in the third plot  only the number of variables fixed for each $s$, considering the  \ref{DDFact} bound, and we can see that the fixing procedure is very effective when $s\leq 80$ .

Our success with fixing using the \ref{DDFact} bound on the $n=2000$ matrix,
for $s=20,40,60,80$, gave us some hope to solve these instances to optimality, or at least reduce them to a size where we could realistically hope that branch-and-bound could succeed. So we devised an iterative
fixing scheme, applying fixing to a sequence of reduced instances, with the goal of solving to optimality
or at least fixing substantially more variables.
At each iteration, we calculated and attempted to fix
based on the \ref{DDFact} bound \emph{and} the \ref{linx} bound. We re-applied the heuristic for a
reduced problem, in
case it could improve on the lower bound of its parent.
Even though the \ref{linx} bound cannot fix any variables at the first iteration, for  $s=20,40,60$,
it enabled us to fix more variables for reduced problems. For $s=20,40,60$, we could solve to
optimality. The results are summarized in Table \ref{table:it_fix}, where $s'$ and $n'$ are the parameters
for reduced problems, by iteration, and
$*$ indicates the iteration
where we can assert that fixing identified the optimal solution. Unfortunately, for $s=80$, \ref{linx} could
not fix anything after one round of \ref{DDFact} fixing.
We noted that for $s=20,40,60$, the heuristic applied to the $n=2000$ root instance gave what turned out to be the optimal solution. It is possible that
we did not succeed on  $s=80$ because our
lower bound is not strong enough.

 We carried out some additional experiments for
the $n=2000$ matrix, with $s=860,880,\ldots,940$ (recall that the rank of
$C$ is 949). For these instances, \ref{DDFact} is very hard
for \texttt{Knitro} to solve: the solution times for \ref{DDFact}
are an order of magnitude larger as compared to the
instances with $20 \leq s\leq 200$; this
is not the case for \ref{linx}. Additionally,
\texttt{Knitro} failed to converge for $s=920$.
In any case, for all of these problems that we
could solve, we had huge integrality gaps,
and no variables could be fixed based on
either \ref{linx} or \ref{DDFact}.

\subsection{More specifics about the computational time}

In Table \ref{tab:wallclocktime}, we show  means and standard deviations of the average wall-clock times (in seconds) for the main procedures considered in our experiments  and for each $n$. For $n=63,90,124$, the statistics consider the solution for all $2\leq s\leq n-1$. For $n=2000$, the statistics consider all $20\leq s\leq 200$ that are multiples of 20.

In Table \ref{tab:wallclocktime}, the columns ``DDFact'', ``DDFactcomp''
and ``linx (Newton)'' summarize information depicted in the second plots of Figures \ref{fig:data63}, \ref{fig:data90}, and \ref{fig:data124}.
For each $n$, the mean and standard deviation of the times  increase in the order: \ref{linx}, \ref{DDFact}, complementary \ref{DDFact},  and the means and standard deviations
are all significantly better for \ref{linx}.

In Table \ref{tab:wallclocktime}, the column ``linx (BFGS)'' presents the mean and standard deviation of the times (across all $s$)
when \ref{linx} is solved by \texttt{Knitro} without passing the Hessian of the objective function to the solver, i.e., with the application of a BFGS-based algorithm. We see that not passing the Hessian of the objective function to the solver leads to a significant increase in the solution time, and also in the variability of the times across the different values of $s$.
Although we can see performance aggregated over $s$ in the ``linx (Newton)'' and ``linx (BFGS)'' columns, in Figure \ref{fig:Linx:time_ratio}, we get a more complete view of the strong dominance, across most $s$ for each input matrix $C$. We have plotted the time for linx (Newton) divided by
time for linx (BFGS), against $s/n$. With the vast majority of the ratios being less than one for each
input matrix, and this emphatically being the case for the $n=2000$ matrix, we can confidently
recommend passing the Hessian to \texttt{Knitro} when solving \ref{linx}.

In Table \ref{tab:wallclocktime}, the column ``$\gamma$ (Newton)'' in Table \ref{tab:wallclocktime}
 presents statistics for the time used to compute the value of the  parameter $\gamma$ used in \ref{linx} across the different $s$.   To optimize $\gamma$ we do a one-dimensional search, exploiting the fact that the \ref{linx} bound is convex  in the logarithm of  $\gamma$ (see  \cite*{Mixing}), and we use \texttt{Knitro} passing the Hessian at each iteration of the one-dimensional search.
We observe that, compared to solving \ref{linx},  optimizing $\gamma$ is very expensive, however, we should note that the optimization procedure applied had no concern with time. When time is relevant, as in the context of a branch-and-bound algorithm, we can apply a faster procedure, like the one applied in \cite*{Kurt_linx}. Furthermore, we should notice  that in a branch-and-bound context, the \ref{linx} bound  is computed for each subproblem considered, but the  parameter $\gamma$ should not be optimized for every one of them. As done in \cite*{Kurt_linx}, it would be more efficient to use the same parameter value as the one used on the parent node most of the times.

\subsection{Some experiments with \ref{CMESP}}
To illustrate the application of our bounds to \ref{CMESP}, we repeated the experiments performed with the instances of \ref{MESP} with covariance matrix of dimension $n=63$ and $5\leq s\leq 47$, but now including five side constraints $a_i^\top x \leq b_i$, for $i=1,\ldots,5$. The left-hand side of constraint $i$
is given by a uniformly-distributed random vector $a_i$ with integer components between $1$ and $5$.
The right-hand side of the constraints was selected so that, for every $5\leq s \leq 47$, the best known solution $x^*(s)$ of the instance of \ref{MESP} is violated by at least one constraint. For that,  each $b_i$ was selected as the  $80$-th percentile of the values $a_i^\top x^*(s)-1$, for all  $5\leq s \leq 47$.
We note that when considering side constraints,  the linear program  \ref{mingapprobmix} does not have a closed-form solution. In this case,  we solve it with \texttt{Knitro}. The time needed to calculate the dual solution with the \texttt{Knitro} linear-programming solver is no more than $5\%$ of the time needed to calculate the \ref{DDFact} bound or the  complementary \ref{DDFact} bound. However, for the \ref{linx} bound, the variability of the times is large;  for some instances, the time needed to construct the dual solution can even exceed the time needed to calculate the \ref{linx} bound.

\begin{figure}[htbp]
	\centering

	\scalebox{0.55}{
		\includegraphics{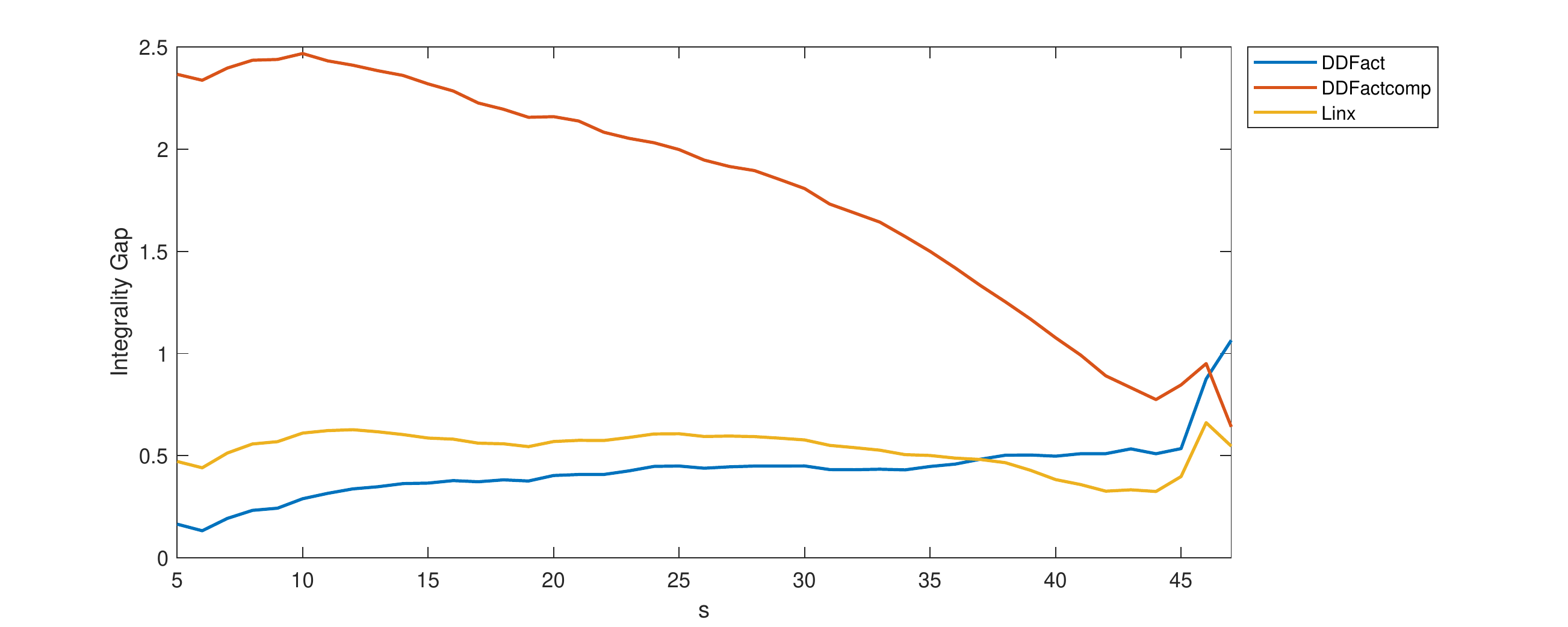}}
	\scalebox{0.55}{
	\includegraphics{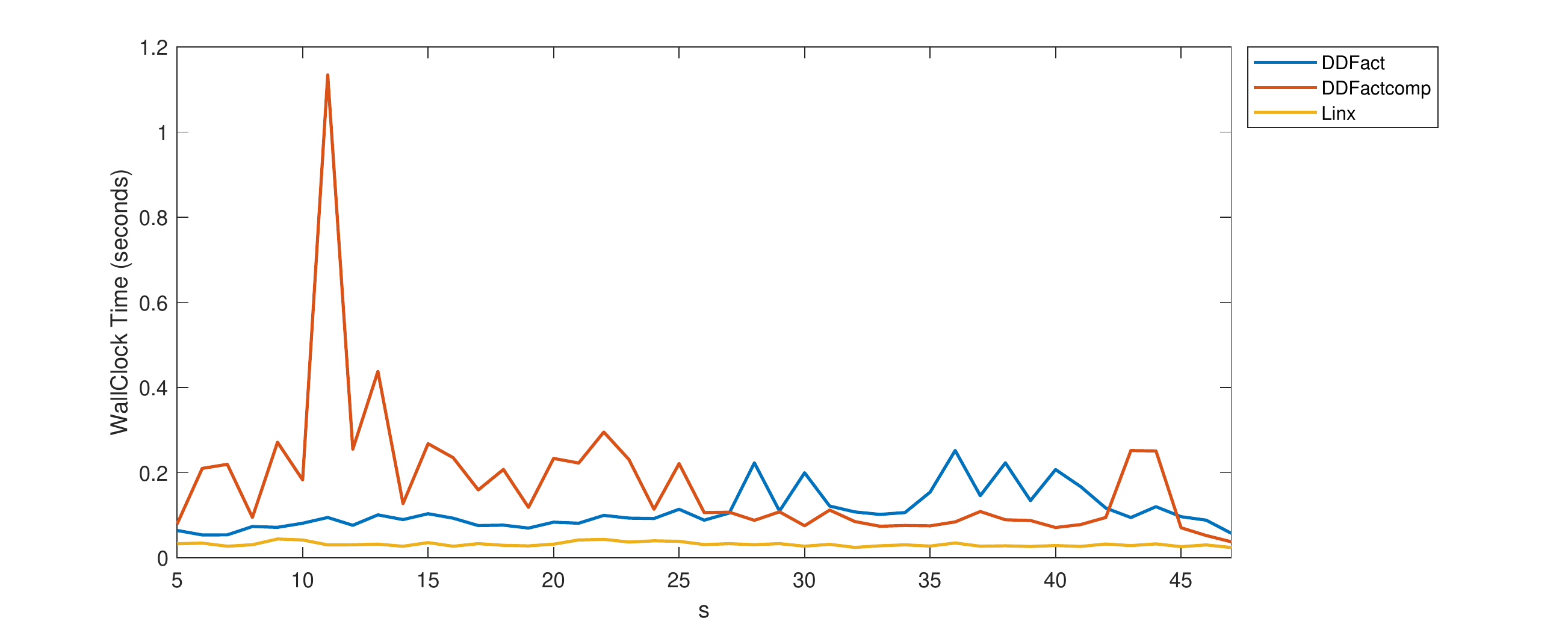}
	}
	\scalebox{0.55}{
	\includegraphics{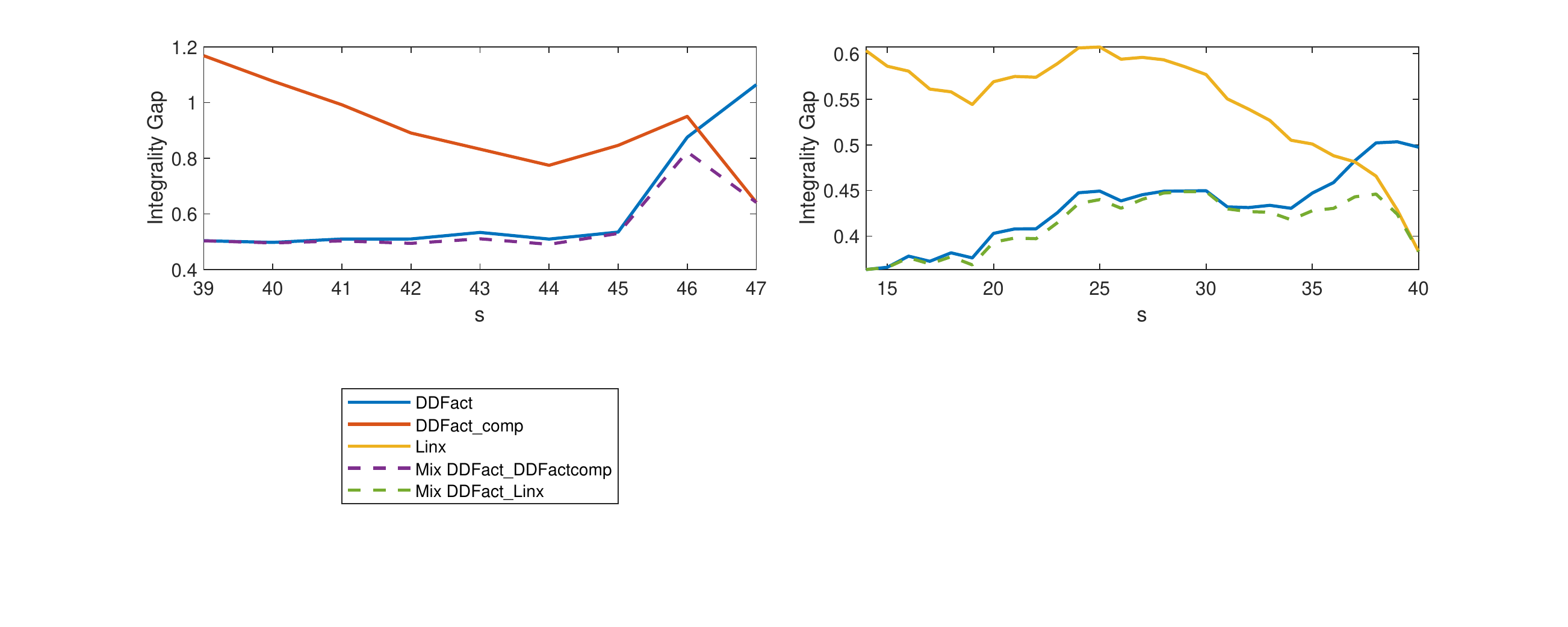}}
	\scalebox{0.55}{
	\includegraphics{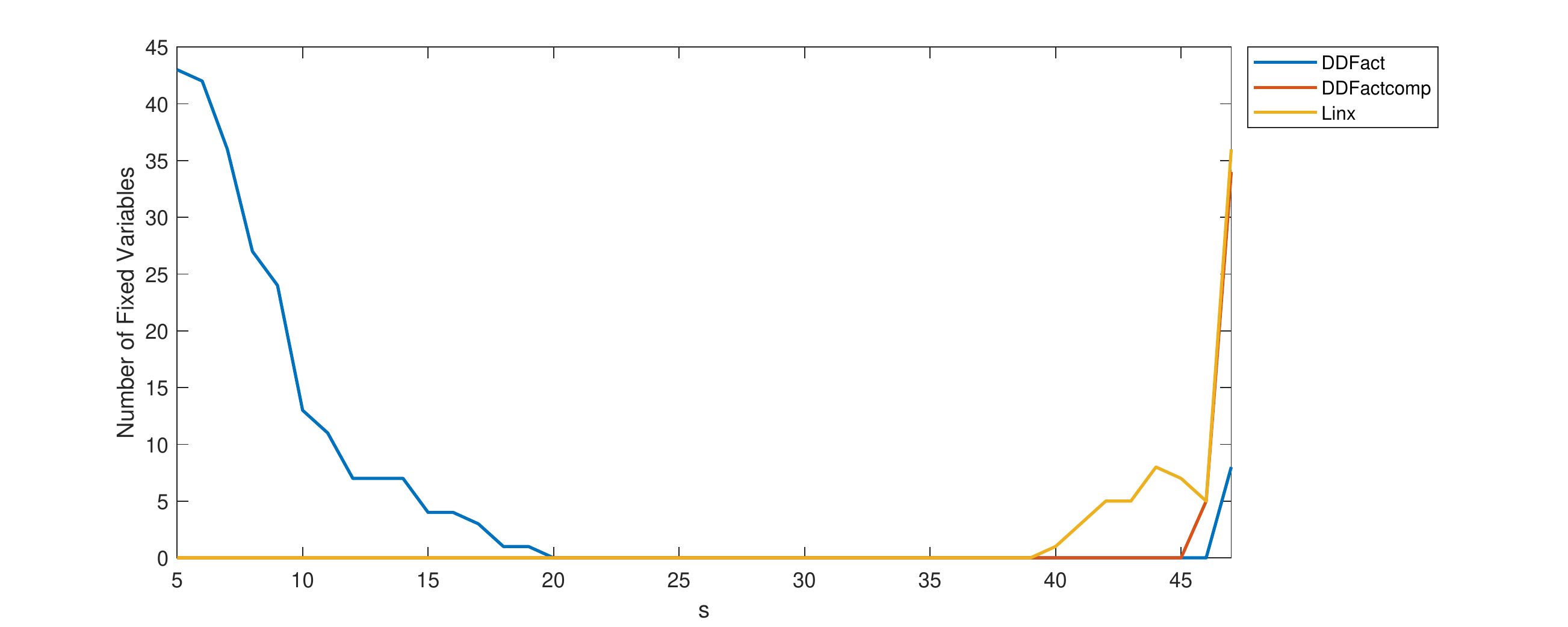}
	}
	\caption{Bounds/times comparison and effect of the mixing  and  variable-fixing methodologies for $n=63$ with 5 side constraints (CMESP)}
	\label{fig:data63_lincon}
\end{figure}

In Figure   \ref{fig:data63_lincon} we present plots for \ref{CMESP}, analogous to those shown in Figure \ref{fig:data63} for \ref{MESP}, considering our instances of dimension $n=63$. We have a very similar analysis of the results presented in both figures, illustrating the robustness of our approach when including  side constraints to \ref{MESP}.


\section{Discussion}\label{sec:discussion}
We developed some useful properties of the
 \ref{DDFact} bound, aimed at guiding computational practice.
 In particular, we saw that (i) the
 \ref{DDFact} bound is invariant under the chosen factorization of the input matrix $C$, (ii)
 the \ref{DDFact} bound cannot be improved by scaling $C$,
 and (iii)  the \ref{DDFact} bound
 dominates the spectral bound.
 We developed a fixing scheme for \ref{DDFact},
 we demonstrated how to mix \ref{DDFact} with
\ref{linx} and with complementary \ref{DDFact} (see \cite*[Sec. 7]{Mixing}) for
general comments on how and when mixing could potentially be employed within branch-and-bound),
and we gave a general variable-fixing scheme for these mixings.

Overall, we found that working with a general-purpose nonlinear-programming solver is quite practical for
solving \ref{linx} and \ref{DDFact}
 relaxations of \ref{CMESP}.
For \ref{DDFact}, this is despite the fact that its objective function is not guaranteed to be smooth at
all iterates (of the nonlinear-programming solver). We found that for \ref{linx}, which has a smooth objective function, passing the Hessian to the
nonlinear-programming solver is quite effective; the running times are much better, in mean and variance, compared to a BFGS-based approach.
We found that various mixings of
\ref{linx}, \ref{DDFact} and complementary \ref{DDFact} can lead to improved bounds. We found that fixing can be quite effective for \ref{DDFact} and complementary \ref{DDFact}.
Unfortunately, we did not find
mixing to be useful for fixing additional variables, as
compared to fixing variables based on each relaxation separately,
on the benchmark instances that we experimented with. But we did find  iterative fixing,
employing the \ref{linx} and \ref{DDFact} fixing rules in concert,
to be quite effective on
 large and difficult instances; for the $n=2000$ matrix and $s=20,40,60$,
we could find and verify optimal solutions for the first time, and without any branching.

In future work, we
plan to develop a full  branch-and-bound implementation
aimed at solving difficult instances of \ref{CMESP} to optimality.
Our experiments on benchmark covariance matrices indicate that such an approach should use both
the \ref{linx} and \ref{DDFact} bounds. In particular \ref{linx}
often seems to be valuable for the large values of $s$ (where \ref{DDFact} deteriorates). We can hope that variable fixing can be exploited for subproblems,
and we expect that mixing is likely to be most valuable for
subproblems that have $s$ near the middle of the range.

Additionally, we plan on working further on improving the algorithmics for the
\ref{DDFact} relaxation, for modern settings in which  the order $n$
of the covariance matrix  greatly exceeds
its rank $r$. Specifically: (i)  we plan to take  better advantage of the
fact that the matrix in the objective function of
\ref{DDFact} is order $r$ while in \ref{linx} it is order $n$; (ii) while the objective
function of \ref{DDFact} is not guaranteed to be smooth at all iterates,
we found that it usually is, and so we plan to use second-order information
to improve convergence.


\newpage

\newpage

\begin{table}[!ht]
    \centering
    \begin{tabular}{C{1.5cm}|C{1.2cm}C{1.2cm}|C{1.2cm}C{1.2cm}|C{1.2cm}C{1.2cm}|C{1.2cm}C{1.2cm}}
\toprule
  \multirow{1}{*}{Iter}
  & $s'$ & $n'$ & $s'$ & $n'$ & $s'$ & $n'$ & $s'$ & $n'$ \\
\cmidrule(lr){1-1}\cmidrule(lr){2-2} \cmidrule(lr){3-3} \cmidrule(lr){4-4} \cmidrule(lr){5-5} \cmidrule(lr){6-6} \cmidrule(lr){7-7}\cmidrule(lr){8-8}\cmidrule(lr){9-9}
    0       &  20  &  2000 &   40 &    2000&  60 &   2000&  80&    2000 \\
    1      &  20  &  28 &   40  &  58 &60  &  110 &  80  & 442\\
    2      & 2  &  7 &   6  &  13 & 22&  68 &   &   \\
    3       &  $*$ & $*$ & 3& 4& 16 & 24& &  \\
    4 &&&1&2&1&2&&\\
    5 &&&$*$&$*$&$*$&$*$&&\\
\bottomrule
\end{tabular}
    \caption{Iterated fixing for $n=2000$}
     \label{table:it_fix}
\end{table}

\FloatBarrier


\begin{table}[ht]
    \centering
    \begin{tabular}{C{1.1cm}|C{1.2cm}C{1.2cm}|C{1.2cm}C{1.2cm}|C{1.2cm}C{1.2cm}|C{1.2cm}C{1.2cm}|C{1.2cm}C{1.2cm}}
\toprule
  \multirow{2}{*}{$n$}& \multicolumn{2}{c|}{DDFact}& \multicolumn{2}{c|}{DDFactcomp}& \multicolumn{2}{c|}{linx (Newton)}& \multicolumn{2}{c|}{linx (BFGS)} &\multicolumn{2}{c}{$\gamma$ (Newton)} \\
 \cmidrule(lr){2-3}\cmidrule(lr){4-5}\cmidrule(lr){6-7} \cmidrule(lr){8-9} \cmidrule(lr){10-11}
  & mean & std & mean & std & mean & std & mean & std & mean & std\\
\cmidrule(lr){1-1}\cmidrule(lr){2-2} \cmidrule(lr){3-3} \cmidrule(lr){4-4} \cmidrule(lr){5-5} \cmidrule(lr){6-6} \cmidrule(lr){7-7}\cmidrule(lr){8-8}\cmidrule(lr){9-9} \cmidrule(lr){10-10} \cmidrule(lr){11-11}
    63       &  0.1807  &  0.0828 &   0.2756 &    0.2496&  0.0459 &   0.0061&  0.0518&    0.0092& 0.4202&  0.0832 \\
    90       &  0.3653  &  0.1327 &   0.4163  &  0.2398 &0.0629  &  0.0094 &  0.0849  & 0.0202&   0.6727  &0.1376\\
    124       & 0.5023  &  0.2373 &   0.7451  &  0.3892 & 0.0874&  0.0142 &  0.1005 &   0.0177 & 1.3043  & 0.2675\\
    2000       &  461.49 & 536.37& -& -& 133.69 & 17.54&542.93 & 530.50 & 2242.90&645.37\\
\bottomrule
\end{tabular}
    \caption{Wallclock time (sec)}
    \label{tab:wallclocktime}
\end{table}

\FloatBarrier

\vbox{
\ACKNOWLEDGMENT{%
J. Lee was supported in part by AFOSR grant FA9550-19-1-0175.
M. Fampa was supported in part by CNPq grants 305444/2019-0 and 434683/2018-3.
}
}

%
\begin{APPENDIX}{Lagrangian dual for mixing}

The formulation \ref{prob:mixing} can be recast as
\begin{equation*}
\max\left\{ \sum_{i=1}^m \alpha_i  f_i(W_i)
 ~:~ \mathbf{e}^{\top}x=s,~ Ax\leq b,~0\leq x\leq \mathbf{e},~ \alpha_i L_i(x)=\alpha_i W_i,~i=1,\ldots,m~ \right\}.
\end{equation*}
Next, we consider the Lagrangian function
\begin{align*}
&\mathcal{L}(W_1,\ldots,W_m,x,\Theta_1,\ldots,\Theta_m,\upsilon,\nu,\pi,\tau):=\\
&\quad \displaystyle \sum_{i=1}^m \alpha_i\Big( f_i(W_i) +\Theta_i\bullet (L_i(x)-W_i)\Big) + \upsilon^\top x +  \nu^\top (\mathbf{e}-x) +\pi^\top (b-Ax) +\tau(s-\mathbf{e}^\top x),
\end{align*}
with $\mbox{dom }\mathcal{L}= \left(\mathbb{S}^{k_i}_{+}\right)^m \!\times\! \mathbb{R}^{n}\!\times\! \left(\mathbb{S}^{k_i}\right)^m\!\times\! \mathbb{R}^{n}\!\times\!  \mathbb{R}^{n}\!\times\!  \mathbb{R}^{m}\!\times\!  \mathbb{R}$. The corresponding dual function is
\begin{equation*}
\mathcal{L}^*(\Theta_1,\ldots,\Theta_m,\upsilon,\nu,\pi,\tau):=\sup_{W_i\succeq 0,~x}\mathcal{L}(W_1,\ldots,W_m,x,\Theta_1,\ldots,\Theta_m,\upsilon,\nu,\pi,\tau),
\end{equation*}
and the Lagrangian dual problem  is
\[
\min\{\mathcal{L}^*(\Theta_1,\ldots,\Theta_m,\upsilon,\nu,\pi,\tau)~:~\upsilon\geq 0,~\nu\geq 0,~\pi\geq 0\}.
\]
We note that
\begin{align*}
&\sup_{W_i\succeq 0,~x}\mathcal{L}(W_1,\ldots,W_m,x,\Theta_1,\ldots,\Theta_m,\upsilon,\nu,\pi,\tau) = \nonumber\\
&\qquad  \displaystyle \sum_{i=1}^m \alpha_i\sup_{W_i\succeq 0}\Big(  f_i(W_i)-\Theta_i\bullet W_i\Big)\\ 
 &\qquad\quad +
 \sup_{x}
 \left\{\displaystyle \sum_{i=1}^m  \alpha_i\Big( \Theta_i\bullet L_i(x)\Big) + \upsilon^\top x +  \nu^\top (\mathbf{e}-x) +\pi^\top (b-Ax) +\tau(s-\mathbf{e}^\top x) \right\}.
\end{align*}
For $i=1,\ldots,m$,   we assume  that for any given $\hat\Theta_i\in\mathbb{S}^{k_i}_{++}$, there is a closed-form solution  $\hat{W_i}$ to  $\sup\{ f_i(W_i)-\hat\Theta_i\bullet W_i~:~W_i\succeq 0\}$, such that $\hat\Theta_i\bullet \hat{W}_i=:\rho_i\in\mathbb{R}$ and   $\Omega_i:\mathbb{S}^{k_i}_{++}\rightarrow\mathbb{R}$, is defined by $\Omega_i(\hat\Theta_i)
:= f_i(\hat{W}_i)$ and is convex. Furthermore, we assume that the supremum is $+\infty$ if $\Theta_i\nsucc 0$. Therefore
\begin{align*}
& \displaystyle \alpha_i \sup_{W_i\succeq 0}\Big(   f_i(W_i)-\Theta_i\bullet W_i\Big) = \left\{
\begin{array}{ll}
\alpha_i\Big( \Omega_i(\Theta_i) - \rho_i\Big) ,&\mbox{ if } \Theta_i\succ 0;\\
+\infty, & \mbox{ otherwise.}
\end{array}
\right.
\end{align*}
We have
\begin{align*}
 &
 \sup_{x}
 \left\{\displaystyle \sum_{i=1}^m  \alpha_i \Big(\Theta_i\bullet L_i(x)\Big) + \upsilon^\top x +  \nu^\top (\mathbf{e}-x) +\pi^\top (b-Ax) +\tau(s-\mathbf{e}^\top x) \right\} \\
&\quad\quad = \left\{\begin{array}{l}
\displaystyle\sum_{i=1}^m \alpha_i\Big(\Theta_i\bullet L_{i0 }\Big) + \nu^\top \mathbf{e}  + \pi^\top b +\tau s,\\
\quad \quad \mbox{if} ~   \displaystyle \sum_{i=1}^m \alpha_i\Big(\Theta_i\bullet L_{ij }\Big) + \upsilon_j  - \nu_j   - \pi^\top A_{\cdot j} - \tau=0,~\mbox{ for } j \in N;\\
+\infty, \quad  \mbox{otherwise.}
\end{array}\right.
\end{align*}
We see that the Lagrangian dual problem  is equivalent to \ref{Dmix}.

\end{APPENDIX}
%
%



\bibliographystyle{informs2014} 
\bibliography{fact_paper} 


\end{document}